\documentclass{article}
\usepackage[left=2cm,top=2cm,right=3cm,bottom=1cm, footskip=.5cm]{geometry}
\usepackage{amsmath}
\usepackage{graphicx}
\usepackage{tikz}
\usepackage{fancyvrb}
\usepackage{pdfpages}
\usepackage[makeroom]{cancel}
\usepackage{amssymb} 
\usepackage[parfill]{parskip}
\usepackage{mathtools}
\usepackage{hyperref}
\usepackage{url}
\usepackage{tikz-cd}
\usepackage{setspace}
\usepackage{amsthm}
\usepackage{xfrac}
\usepackage{color}
\usepackage{cite}
\usetikzlibrary{%
	matrix,%
	calc,%
	arrows%
}

\newtheoremstyle{colon}%
{}
{}
{\itshape}
{}
{\bfseries}
{:}
{ }
{}

\theoremstyle{colon}

\newtheorem{Lemma}{Lemma}[section]
\newtheorem{Theorem}[Lemma]{Theorem}

\newtheorem{Corollary}[Lemma]{Corollary}
\newtheorem{Definition}[Lemma]{Definition}
\newtheorem{Example}[Lemma]{Example}
\newtheorem{Remark}[Lemma]{Remark}

\title{Localisations and completions of nilpotent G-spaces}

\begin{document}
	
	\date{}
	\author{A. Ronan}
	\maketitle
	
	\newcommand{\orb}[2][G]{[\sfrac{#1}{#2}]}

\begin{abstract}
	We develop the theory of nilpotent $G$-spaces and their localisations, for $G$ a compact Lie group, via reduction to the non-equivariant case using Bousfield localisation. One point of interest in the equivariant setting is that we can choose to localise or complete at different sets of primes at different fixed point spaces - and the theory works out just as well provided that you invert more primes at $K \leq G$ than at $H \leq G$, whenever $K$ is subconjugate to $H$ in $G$. We also develop the theory in an unbased context, allowing us to extend the theory to $G$-spaces which are not $G$-connected.
\end{abstract} 

\section{Introduction}

The purpose of this paper is to develop the theory of localisations and completions of nilpotent $G$-spaces at sets of primes, where $G$ is a compact Lie group. The main reference for the equivariant theory is \cite[Ch. II]{M96}, which itself is a summary of the older papers \cite{MMT82} and \cite{M82}, where it was explained how the foundations of the theory could be developed using the same arguments as in the non-equivariant setting, with some additional complications when $G$ is compact Lie rather than just finite. Our approach is slightly different, in that we use the theory of Bousfield localisation to deduce the foundations of the theory from the non-equivariant case. This approach leads to fewer difficulties in the compact Lie case, and allows us to use a more general definition of a nilpotent $G$-space than in \cite{M96}, see Definition \ref{name21}. For example, we prove that a nilpotent $G$-space $X$ is $p$-complete iff all homotopy groups of the form $\pi_i(X^H)$ are $p$-complete. This fact was proved in \cite[Theorem 2]{M82}, but only under the assumption that, for fixed $i$, the nilpotency classes of $\pi_i(X^H)$, as $H$ varies, have a common bound. 

Another contribution of this paper is that we allow the set of primes we are localising or completing at to vary over the orbit category of $G$. In this case, the theory works out just as well provided that you `invert more primes' at $K \leq G$ than at $H \leq G$ whenever $K$ is subconjugate to $H$ in $G$. For example, we could localise at $p$ at one subgroup and complete at $p$ at another, where, loosely speaking, completing at $p$ `inverts more primes' than localising at $p$. One might ask, why consider these localisations? In the non-equivariant setting, Bousfield proved in \cite[Theorem 1.1]{B74} that all localisations at connective homology theories are equivalent to localisations with respect to either $H(-;\mathbb{Z}_T)$ or $H(-; \oplus_{p \in T} \mathbb{F}_p)$ for some set of primes $T$. Therefore, in this paper we are considering localisations at pointwise connective homology theories, where pointwise means we choose a connective homology theory for every closed subgroup $H$ of $G$, and the localisations which satisfy the poset condition are precisely those with the property that a $G$-space is local iff it is pointwise local. 

We develop the theory in both a based and unbased context - with different parts of the theory working better in each setting. For example, we derive fracture theorems for virtually nilpotent $G$-spaces in Theorems \ref{name88} and \ref{name91}, relate localisations of nilpotent $G$-spaces to equivariant Postnikov towers, and show that our homological approach to the theory is equivalent to the classical cohomological approach of \cite{MMT82} and \cite{M82}, all in the based context. We use the unbased theory to extend our results on nilpotent $G$-spaces to $G$-spaces whose fixed point spaces are disjoint unions of nilpotent spaces. This is especially pertinent in the equivariant setting, since there are many examples of $G$-spaces which are non-equivariantly connected, but which have disconnected fixed point spaces, or no possible choice of a $G$-fixed basepoint at all.

We conclude this introduction by describing some of the results from Section \ref{name187} in more detail. Since we intend to derive fracture theorems for virtually nilpotent $G$-spaces, we include a brief but self-contained introduction to the theory of virtually nilpotent spaces in Subsection \ref{name182}. We work with homological localisations throughout, so this introduction should be accessible to all readers, regardless of whether they are familiar with the theory of Bousfield and Kan's $R$-completion functor. Our main fracture theorem for virtually nilpotent $G$-spaces is then as follows:

\begin{Theorem} 
	Consider a commutative square of virtually nilpotent $G$-spaces:
	\begin{center}

		\begin{tikzcd}
			X \arrow{r}{f} \arrow[swap]{d}{\phi} & Y \arrow{d}{\psi} \\
			A \arrow[swap]{r}{g} & B \\
		\end{tikzcd}
	\end{center}	
	
	such that, for every subgroup $H$ of $G$, there are sets of primes $T_H, S_H$ satisfying:
	
	i) $X^H,Y^H$ are $T_H$-local and $A^H,B^H$ are $S_H$-local, \\
	ii) $f^H$ is an $\mathbb{F}_{T_H}$-equivalence and $g^H$ is an $\mathbb{F}_{S_H}$-equivalence, \\
	iii) $\phi^{H}, \psi^{H}$ are $\mathbb{Q}$-equivalences.
	
	Then the square is a homotopy fibre square.	
\end{Theorem}

Note that non-equivariantly this fracture theorem tells us, amongst other things, that the local fracture theorem of \cite[Theorem 8.1.3]{MP12} also continues to hold in the virtually nilpotent context. We move on to derive fracture theorems for homotopy classes in Subsection \ref{name196}, the main results being Theorem \ref{name91} and \ref{name188}. A particularly interesting special case of Theorem \ref{name188} is the following fracture theorem for a class of non-constant localisation systems:

\begin{Theorem}
	Let $T$ be a set of primes and suppose that $\mathbf{T}$ is a localisation system whose underlying set of primes is constantly $T$. Let $X$ be a nilpotent $G$-space such that, for every $H$ and $i \geq 1$, $\pi_i(X^H)$ is $f \mathbb{Z}_T$-nilpotent.  Then, for every finite based $G$-CW complex $K$, there is a pullback of sets:
	
		\[		
	\begin{tikzcd}
		{[K,X]} \arrow{r}{} \arrow[swap]{d}{} & {[K,X_{\mathbf{T}}]} \arrow{d}{} \\
		{[K,X_0]} \arrow[swap]{r}{} & {[K,(X_\mathbf{T})_0]} 
	\end{tikzcd}
	\]
\end{Theorem}
The proof of Theorem \ref{name188} is surprisingly simple, although the proof of a closely related algebraic result, Lemma \ref{name193}, is more involved. We close Subsection \ref{name196} with an analysis of what goes wrong when we allow the underlying set of primes of $\mathbf{T}$ to vary over the orbit category. This leads to a particularly nice application of the theory - namely, the construction of maps $f,g: K \to X$ which are not equivariantly homotopic even though, for all $H$, $f^H$ and $g^H$ are homotopic. Here we can take both $K$ and $X$ to be $\mathbf{T}$-localisations of finite $G$-CW complexes, where the coefficient of $\mathbf{T}$ is constantly $0$ and provided that the underlying set of primes of $\mathbf{T}$ is not constant.

We include a brief discussion of equivariant phantom maps in Subsection \ref{name186}. We develop the theory of equivariant Postnikov towers in Subsection \ref{name80}, obtaining the following sharp description of Postnikov towers in Theorem \ref{name185}:

\begin{Theorem}
	A nilpotent $G$-space is $\mathbf{T}$-local iff it is $\mathcal{B}_{\mathbf{T}}$-nilpotent, where $\mathcal{B}_{\mathbf{T}}$ is the class of $\mathbf{T}$-local coefficient systems. A $G$-space is a $\mathbf{T}$-local bounded $\mathcal{A}$-nilpotent $G$-space iff it is bounded $\mathcal{B}_{\mathbf{T}}$-nilpotent. 
\end{Theorem}

Finally, we conclude the paper by applying the theory of equivariant Postnikov towers to prove the following result, which is Theorem \ref{name29}: 

\begin{Theorem}  A map of $G$-spaces, $f:X \to Y$, is a $\mathbf{T}$-equivalence iff it is a cohomology $\mathbf{T}$-equivalence.	
\end{Theorem}

Therefore, our homological approach to the theory is equivalent to the classical cohomological approach of \cite{MMT82} and \cite{M82}.

\subsection{Notations and Prerequisites}

We will work with the model categories of $G$-spaces and based $G$-spaces, where $G$ is a compact Lie group, basepoints are $G$-fixed, and the model structures are the Quillen or $q$-model structures, \cite[Proposition B.7]{S18}, \cite[Theorem 7.6.5]{H03}. All subgroups of $G$ under consideration are assumed to be closed. Unless otherwise stated, we build $G$-CW complexes out of the maps $(\frac{G}{H})_+ \wedge S^n_+ \to (\frac{G}{H})_+ \wedge D^n_+$ in the based context, rather than using based maps out of $(\frac{G}{H})_+ \wedge S^n$. The basepoint of a based $G$-space, $X$, is said to be nondegenerate if $* \to X$ is a $G$-Hurewicz cofibration (a map satisfying the equivariant analogue of the HEP, \cite[pg. 43]{M99}). The notation $[A,B]$ denotes homotopy classes of maps, which may be based/unbased/equivariant depending on the context.

This paper should be accessible to any reader who is familiar with the non-equivariant theory of nilpotent spaces and their localisations, as well as the basics of equivariant homotopy theory.

\section{Localisation systems}

\subsection{Bousfield localisation at the T-equivalences} \label{name50}

In this subsection, we define localisation systems, $\mathbf{T}$, as well as the notion of a $\mathbf{T}$-equivalence between based $G$-spaces. We develop the basic properties of the $\mathbf{T}$-equivalences, and then use the Bousfield cardinality argument to show that there exists a model structure on the category of based $G$-spaces, where a map is a weak equivalence iff it is a $\mathbf{T}$-equivalence. We then derive some properties of this model structure, including Theorem \ref{name100} below, which is the key to deducing the equivariant theory of nilpotent $G$-spaces from the non-equivariant theory.

Let $\mathbf{P}$ denote the poset of subsets of the set of prime numbers partially ordered by inclusion, and let $\mathcal{O}$ denote the orbit category of a compact Lie group $G$. We begin with the following important definitions: 

\begin{Definition}
	A localisation system is a functor $\mathbf{T}: \mathcal{O}^{op} \to \mathbf{P}^{op} \times \mathbf{1}$, where we denote by $\mathbf{1}$ the category with objects 0 and 1 and a single arrow from 0 to 1.
\end{Definition}

We think of $\mathbf{T}(\orb{H})$ as a set of primes with coefficient, where the coefficient is either 0 or 1. Intuitively, at least for nilpotent spaces, a coefficient of 0 will mean we are localising at $T$, and a coefficient of 1 will mean we are completing at $T$. If we drop the bold font on the $\mathbf{T}$, then $T(\orb{H})$ denotes only the underlying set of primes of $\mathbf{T}(\orb{H})$. More explicitly:

\begin{Definition}
	A set of primes with coefficient is a pair $(T,n)$ where $T$ is a set of primes and $n = 0$ or $1$. We call $T$ the underlying set of primes and $n$ the coefficient. 
\end{Definition}

 Let $\mathbb{Z}_T$ denote the integers localised at $T$, where recall that localising at $T$ inverts the primes not in $T$. Recall that a map of spaces is called a $\mathbb{Z}_T$-equivalence if it induces an isomorphism on homology with coefficients in $\mathbb{Z}_T$. Similarly, a map is called an $\mathbb{F}_T$-equivalence if it induces an isomorphism on homology with coefficients in $\mathbb{F}_p$, for every $p \in T$. When the basepoints are nondegenerate, it is equivalent to define these equivalences using the respective reduced homology theories instead.

\begin{Definition}
	Let $\mathbf{T}$ be a set of primes with coefficient and $f: X \to Y$ a map of spaces. If the coefficient is 0, then we call $f$ a $\mathbf{T}$-equivalence if it is a $\mathbb{Z}_T$-equivalence. If the coefficient is 1, then we call $f$ a $\mathbf{T}$-equivalence if it is an $\mathbb{F}_T$-equivalence.
\end{Definition}

We have the following similar definition in the case where $\mathbf{T}$ is a localisation system. Note that if $G$ is the trivial group, so we are working non-equivariantly, then a localisation system is equivalent to a set of primes with coefficient, and the two definitions are consistent with one another:

\begin{Definition} \label{misc219}
	Let $\mathbf{T}$ be a localisation system and $f: X \to Y$ be a map of based $G$-spaces. We say that $f$ is a $\mathbf{T}$-equivalence if for all $H \leq G$, $f^{H}: X^H \to Y^H$ is a $\mathbf{T}(\orb{H})$-equivalence. 
\end{Definition}

We now define $\mathbf{T}$-local spaces. Since we are working with homological localisations in this section, we should be precise with our terminology in order to avoid confusion with the $R$-completion functor of \cite{BK72}. Therefore, if $T$ is a set of primes, by a $T$-local space we always mean a fibrant object in the Bousfield localised model structure with respect to the homology theory $H(-;\mathbb{Z}_T)$, \cite[Theorem 19.3.8]{MP12}. Moreover:

\begin{Definition}
	Let $\mathbf{T}$ be a set of primes with coefficient and let $X$ be a space. If the coefficient is 0, we say that $X$ is $\mathbf{T}$-local if it is $T$-local. If the coefficient is $1$, we say that $X$ is $\mathbf{T}$-local if it is a fibrant object in the Bousfield localised model structure with respect to the homology theory $H(-;\oplus_{p \in T} \mathbb{F}_p)$.
\end{Definition}

Our next task is to Bousfield localise the category of $G$-spaces with respect to the class of $\mathbf{T}$-equivalences, where $\mathbf{T}$ is a localisation system. However, to aid clarity at certain points of this section, we will consider, at first, the more general situation of Bousfield localisation with respect to the class of $\mathbf{A}$-equivalences, as defined below, where $\mathbf{A}$ is any function from the set of objects, $\orb{H}$, of the orbit category to abelian groups. 

\begin{Definition}
	Let $\mathbf{A}$ be a function from the set of objects of the orbit category of $G$, or equivalently from the set of subgroups of $G$, to abelian groups. A map of $G$-spaces, $f: X \to Y$, is said to be an $\mathbf{A}$-equivalence if, for all subgroups $H$ of $G$, $f^H$ induces an isomorphism on homology with coefficients in $\mathbf{A}(\orb{H})$.
\end{Definition}

Note that, if $\mathbf{T}$ is a localisation system, then a $\mathbf{T}$-equivalence is equivalent to an $\mathbf{A}$-equivalence, for some function $\mathbf{A}$ taking values in abelian groups of the form $\mathbb{Z}_T$ or $\oplus_{p \in T} \mathbb{F}_p$.

We will need the following minimal list of properties of the $\mathbf{A}$-equivalences, where a property is \textit{pointwise} if it holds for all fixed point spaces:

\begin{Lemma} \label{name4}
	i) The class of $\mathbf{A}$-equivalences is closed under retracts, satisfies 2-out-of-3, and every weak equivalence is an $\mathbf{A}$-equivalence, \\
	ii) the pushout of an $\mathbf{A}$-equivalence that is a pointwise $h$-cofibration is an $\mathbf{A}$-equivalence, \\
	iii) the colimit of a transfinite sequence of $\mathbf{A}$-equivalences which are closed inclusions is an $\mathbf{A}$-equivalence. 
\end{Lemma}

\begin{proof}	
	i) is easy. For ii), since taking fixed points preserves pushouts along closed inclusions, \cite[Proposition B.1]{S18}, we can work pointwise and replace the spaces with nondegenerately based ones. The result then follows from consideration of cofibre sequences. For iii), taking fixed points preserves transfinite colimits of closed inclusions, \cite[Proposition B.1]{S18}, and so the result follows from the fact that homology preserves these colimits. 
\end{proof}	

We can use the Bousfield-Smith cardinality argument on the $\mathbf{A}$-equivalences. The argument is essentially the same as the classical case of localising spaces with respect to homology theories, which is treated in \cite[Section 19.3]{MP12}. The key lemma is as follows, where all cell complexes are $G$-cell complexes:

\begin{Lemma} \label{name5}
	There exists a cardinal $\kappa$ with the following property: if $i:A \to B$ is the inclusion of a subcomplex into a cell complex $B$ which is also an $\mathbf{A}$-equivalence, then, for any cell $e$ of $B$, there is a subcomplex $C$ of size $< \kappa$ containing $e$ such that $A \cap C \to C$ is an $\mathbf{A}$-equivalence. 
\end{Lemma}

\begin{proof}
	Choose a regular cardinal, $\kappa > \aleph_0$, with the following properties:
	
	i) every cell of any cell complex is contained within a subcomplex of size  $< \kappa$, \\
	ii) if $Z$ is a cell complex of size $< \kappa$, then $\bigoplus_{*,H} H_*(Z^H; \mathbf{A}(\orb{H}))$ has cardinality $< \kappa$, \\	
	iii) if $W$ is any cell complex, then for any $*$ and $H$: 
	\[ H_*(W^H; \mathbf{A}(\orb{H})) = \mathrm{colim}_{<\kappa} H_*(Z^H; \mathbf{A}(\orb{H}))
	\]
	where the colimit is over all subcomplexes of $W$ of size $< \kappa$.
	
To start the proof, choose a subcomplex $C_0$ of B of size $< \kappa$ which contains $e$, so we have a map $C_0 \cap A \to C_0$. By some $c \in H_*(C_0, C_0 \cap A; \mathbf{A})$, we mean an element of $H_n(C_0^H, (C_0 \cap A)^H; \mathbf{A}(\orb{H}))$ for some $H$ and $n$. For each $c \in H_*(C_0, C_0 \cap A; \mathbf{A})$, its image in $H_*(B,A; \mathbf{A}) = 0$ vanishes. We know that $H_*(B,A; \mathbf{A}) = colim_{< \kappa} H_*(B^{'},B^{'} \cap A; \mathbf{A})$, where the colimit runs over all subcomplexes, $B^{'}$, of $B$ of size $< \kappa$. Therefore, there exists a $< \kappa$ dimensional subcomplex $E$ of $B$, containing $C_0$, such that the image of $c$ in $H_*(E, E \cap A; \mathbf{A})$ is $0$. Define $C_1$ by adding such a subcomplex $E$ to $C_0$ for every $c \in H_*(C_0, C_0 \cap A; \mathbf{A})$ - the conditions i) - iii) above ensure that $C_1$ has size $< \kappa$. Similarly, define $C_2,C_3,...$ and define $C$ to be the union of the $C_i$ which still has size $< \kappa$. Then, since homology preserves such sequential colimits, we have $H_*(C, C \cap A; \mathbf{A}) = 0$, as desired.	
\end{proof}

We now deduce the standard consequences of Lemma \ref{name5}. Firstly, using transfinite induction, we have:

\begin{Corollary}
	A map has the RLP with respect to all inclusions of cell complexes which are $\mathbf{A}$-equivalences iff it has the RLP with respect to all inclusions of cell complexes of dimension $< \kappa$ which are $\mathbf{A}$-equivalences.
\end{Corollary}

\begin{proof}
	See \cite[Proposition 4.5.6]{H03}.
\end{proof}

Any map with the RLP with respect to inclusions of cell complexes that are $\mathbf{A}$-equivalences is a $q$-fibration (see Notations and Prerequisites), since the generating acyclic cofibrations $(\frac{G}{H})_+ \wedge (D^n)_+ \to (\frac{G}{H})_+ \wedge (D^n \times I)_+$ are inclusions of subcomplexes. Therefore, using left properness we have:

\begin{Lemma}\label{name7}
	A map has the RLP with respect to all $q$-cofibrations which are $\mathbf{A}$-equivalences iff it has the RLP property with respect to all inclusions of cell complexes that are $\mathbf{A}$-equivalences.
\end{Lemma}

\begin{proof} See \cite[Proposition 13.2.1]{H03}. 
\end{proof}

If we call such a map an $\mathbf{A}$-fibration, then we see that an $\mathbf{A}$-fibration that is an $\mathbf{A}$-equivalence is a $q$-acyclic $q$-fibration by the retract argument. Using this and the small object argument we can now conclude:

\begin{Theorem} \label{name8}
	There is a left proper model structure on the category of based $G$-spaces where the weak equivalences are the $\mathbf{A}$-equivalences, the cofibrations are the $q$-cofibrations and the fibrations are the $\mathbf{A}$-fibrations.
\end{Theorem}

This model structure is monoidal:

\begin{Lemma} \label{name81} 
	If $i:A \to B$ and $j:C \to D$ are cofibrations, then $i \square j: A \wedge D \cup B \wedge C \to B \wedge D$ is a cofibration which is an $\mathbf{A}$-equivalence if either i or j is an $\mathbf{A}$-equivalence.
\end{Lemma}

\begin{proof}
	The fact that $i \square j$ is a cofibration is classical and is a consequence of the fact that $\frac{G}{H} \times \frac{G}{K}$ is $G$-homeomorphic to a $G$-CW complex. Similarly, since $(\frac{G}{H})^K$ is homeomorphic to a CW-complex, by \cite[Corollary 7.2]{I83} and \cite[Ch. VI, Corollary 2.5]{B72}, we have that a cofibration is a pointwise cofibration. Therefore, for the remaining statement concerning $\mathbf{A}$-equivalences we can assume that $G$ is the trivial group. Note also that the cofibre of $i \square j$ is homotopy equivalent to $\frac{B}{A} \wedge \frac{D}{C}$. Suppose that $j$ is an $\mathbf{A}$-equivalence and $p$ is an $\mathbf{A}$-fibration. Then $i \square j$ has the left lifting property with respect to $p$ iff $i$ has the left lifting property with respect to $p^{\square j}$. Therefore, it suffices to show that $i \square j$ is an $\mathbf{A}$-equivalence in the case where $i: (\frac{G}{H})_+ \wedge (S^{n-1})_+ \to (\frac{G}{H})_+ \wedge (D^n)_+$ and $j$ is an inclusion of a subcomplex which is an $\mathbf{A}$-equivalence. Since we are assuming that $G$ is trivial, from the third sentence of this proof, the cofibre of $i \square j$ is homotopy equivalent to $S^n \wedge \frac{D}{C}$, which has vanishing reduced homology with the required coefficients as desired. 
\end{proof}

We have the following characterisation of the fibrant objects:

\begin{Lemma} \label{name56}
	Let $Z$ be a based $G$-space. Then the following are equivalent:
	
	i) $Z$ is $\mathbf{A}$-local (that is fibrant in the model structure of Theorem \ref{name8}), \\	
	ii) for all $\mathbf{A}$-equivalences $f: A \to B$ between cofibrant objects, the map $f^*: Map(B,Z) \to Map(A,Z)$ is a $G$-weak equivalence, \\
	iii) for all $\mathbf{A}$-equivalences $f: A \to B$ between cofibrant objects, the map $f^*: [B,Z] \to [A,Z]$ is a bijection.
\end{Lemma}

\begin{proof}	
	$i) \implies ii)$ We can assume that $f$ is a cofibration. It suffices to show that $f^*$ is a $q$-acyclic $q$-fibration, which follows from Lemma \ref{name81} since a lifting problem against $f^*$ transposes to a lifting problem against $Z \to *$, as in \cite[Proposition 12.2]{B75}. \\
	$ii) \implies iii)$ follows from the fact that $Map(A,Z)^G = Map_G(A,Z)$ and passage to $\pi_0$, \\
	$iii) \implies i)$ it is easy to see that $Z \to  *$ has the right lifting property with respect to any inclusion of cell complexes that is an $\mathbf{A}$-equivalence, using the fact that inclusions of cell complexes are $h$-cofibrations, and so $Z$ is $\mathbf{A}$-local by Lemma \ref{name7}.
\end{proof}

We record the following closure properties of $\mathbf{A}$-local $G$-spaces since we will use their non-equivariant analogues, \cite[Theorem 12.9]{B75}, in Subsection \ref{name182}:

\begin{Lemma} \label{name183}
	i) If $f: X \to Z$ and $g: Y \to Z$ are maps between $\mathbf{A}$-local $G$-spaces, then the double mapping path space $N(f,g)$ is $\mathbf{A}$-local, \\
	ii) If $f_i : A_{i+1} \to A_i$ is a sequence of maps between $\mathbf{A}$-local $G$-spaces, then $holim A_i$ is $\mathbf{A}$-local.
\end{Lemma}

\begin{proof}
	Both are consequences of the fact that a $q$-fibration between $\mathbf{A}$-local objects is an $\mathbf{A}$-fibration, \cite[Proposition 3.3.16]{H03}.
\end{proof}

From now on, we return to the context of Bousfield localisation at the $\mathbf{T}$-equivalences, where $\mathbf{T}$ is a localisation system. In this context, by a $\mathbf{T}$-fibration etc., we mean an $\mathbf{A}$-fibration, as above, for the function $\mathbf{A}$ defined by $\mathbf{T}$. The next lemma is the key to deducing our results on nilpotent $G$-spaces from the non-equivariant theory, and is also the first to make use of the definition of a localisation system:

\begin{Lemma} \label{name11}
	Let $\mathbf{T}$ be a localisation system. If a based $G$-space $Z$ is $\mathbf{T}$-local, then $Z^H$ is $\mathbf{T}(\orb{H})$-local for every $H \leq G$. 
\end{Lemma} 

\begin{proof}
	Let $f: A \to B$ be a $\mathbf{T}(\orb{H})$-equivalence between cofibrant spaces. Note that $A$ and $B$ are just spaces, not $G$-spaces. Let $g = 1 \wedge f: (\sfrac{G}{H})_+ \wedge A \to (\sfrac{G}{H})_+ \wedge B$. We have $(\sfrac{G}{H})^K_+ = \mathcal{O}(\orb{K}, \orb{H})_+$, and it follows that $g$ is a $\mathbf{T}$-equivalence by Lemma \ref{name81}, the fact that $\mathbf{T}$ is a localisation system and the following observations:
	
	i) if $S \subset T$, then a $\mathbb{Z}_T$-equivalence is a $\mathbb{Z}_S$-equivalence, \\
	ii) if $S \subset T$, then an $\mathbb{F}_T$-equivalence is an $\mathbb{F}_S$-equivalence, \\
	iii) a $\mathbb{Z}_T$-equivalence is an $\mathbb{F}_T$-equivalence. 
	
	Moreover, since $Z$ is $\mathbf{T}$-local and $g$ is a $\mathbf{T}$-equivalence, we have that $g^*: [(\sfrac{G}{H})_+ \wedge B,Z] \to [(\sfrac{G}{H})_+ \wedge A,Z]$ is a bijection. This is equivalent to $[B,Z^H] \to [A, Z^H]$ being a bijection, and it follows that $Z^H$ is $\mathbf{T}(\orb{H})$-local. 	
\end{proof}

Using Lemma \ref{name11}, we can deduce the following lemma, where by a $\mathbf{T}$-localisation we mean a fibrant replacement in the model structure of Theorem \ref{name8}:

\begin{Theorem} \label{name100}
	i) A based $G$-space $Z$ is $\mathbf{T}$-local iff $Z^H$ is $\mathbf{T}(\orb{H})$-local for every $H \leq G$, \\
	ii) A map of based $G$-spaces $X \to Y$ is a $\mathbf{T}$-localisation iff $X^H \to Y^H$ is a $\mathbf{T}(\orb{H})$-localisation for every $H \leq G$.
	
\end{Theorem}

\begin{proof}
	i) For the direction we haven't already proved, let $Z$ be a based $G$-space such that $Z^H$ is $\mathbf{T}(\orb{H})$-local for every $H \leq G$. Consider a $\mathbf{T}$-localisation $Z \to W$. Then, each map $Z^H \to W^H$ is a $\mathbf{T}(\orb{H})$-equivalence between $\mathbf{T}(\orb{H})$-local objects, and so a weak equivalence as desired. \\
	ii) This follows from i). 
\end{proof}

To end this subsection, we quickly give a counterexample to indicate what can happen if $\mathbf{T}$ is not a localisation system.

\begin{Example}
 Let $G = C_2$, and define a $\mathbf{T}$-equivalence to be a map of based $G$-spaces, $f: X \to Y$, such that $H_*(f^{e}; \mathbb{Z}[p^{-1}])$ and $H_*(f^{G}; \mathbb{Z}[p^{-1},q^{-1}])$ are isomorphisms, where $p$ and $q$ are distinct primes. If $\mathbf{T}$-local $G$-spaces were always pointwise local, then the analogue of Theorem \ref{name100} would also have to hold. Let $\underline{\mathbb{Z}[p^{-1}]}$ denote the constant coefficient system and consider the map $K(\underline{\mathbb{Z}[p^{-1}]},1) \to K(\underline{\mathbb{Z}[p^{-1}]},1)_{\mathbf{T}}$. The induced map between systems of homotopy groups would result in a commutative triangle:

\[
\begin{tikzcd}
	\mathbb{Z}[p^{-1}] \arrow[swap]{d}{\phi} \arrow{dr}{1} & \\
	\mathbb{Z}[p^{-1},q^{-1}] \arrow{r} & \mathbb{Z}[p^{-1}] \\
\end{tikzcd}
\]

where $\phi$ denotes localisation, which is a contradiction since the bottom map has to be the zero map. Therefore, a $\mathbf{T}$-local space is not necessarily pointwise local. \end{Example}

\subsection{Unbased $T$-localisations}

The theory described in Subsection \ref{name50} goes through essentially unchanged in the unbased context. Recall that we are using the same notation for based and unbased homotopy classes of maps. We have:

\begin{Theorem} \label{name51}
	Let $\mathbf{T}$ be a localisation system. There is a left proper, monoidal model structure on the category of $G$-spaces where the weak equivalences are the $\mathbf{T}$-equivalences, the cofibrations are the $q$-cofibrations and the fibrations are the $\mathbf{T}$-fibrations (which are defined as in Lemma \ref{name7}). A $G$-space $Z$ is $\mathbf{T}$-local (that is fibrant in this model structure) iff for all $\mathbf{T}$-equivalences $f: A \to B$ between cofibrant objects, the map $f^*: [B,Z] \to [A,Z]$ is a bijection.
\end{Theorem}

\begin{proof} The existence of the left proper model structure follows from the Bousfield cardinality argument, as in Subsection \ref{name50}. If $i:A \to B$ is a cofibration and $f: X \to Y$ is a cofibration which is a $\mathbf{T}$-equivalence, then $A \times Y \cup B \times X \to B \times Y$ is a cofibration as in Lemma \ref{name81} and it will be a $\mathbf{T}$-equivalence if $(A \times Y \cup B \times X)_+ \to (B \times Y)_+$ is a $\mathbf{T}$-equivalence. The latter map can be identified with $i_+ \square f_+$, which is a $\mathbf{T}$-equivalence by Lemma \ref{name81}. The characterisation of the fibrant objects now follows as in Lemma \ref{name56}. \end{proof}

We define a $\mathbf{T}$-localisation to be a fibrant replacement in the model structure of Theorem \ref{name51}. Since $\mathbf{T}$ is a localisation system, the arguments of Lemma \ref{name11} and Theorem \ref{name100} show:

\begin{Theorem} \label{name53}
	i) A $G$-space $Z$ is $\mathbf{T}$-local iff $Z^H$ is $\mathbf{T}(\orb{H})$-local for every $H \leq G$, \\
	ii) A map of $G$-spaces $X \to Y$ is a $\mathbf{T}$-localisation iff $X^H \to Y^H$ is a $\mathbf{T}(\orb{H})$-localisation for every $H \leq G$.	
\end{Theorem}

At this point, it is helpful to compare based and unbased localisations in the non-equivariant setting. In this setting, a localisation system, $\mathbf{T}$, is equivalent to a set of primes with coefficient. We have:

\begin{Lemma} \label{name54} Let $Z$ be an unbased space. Then: \\
	i) $Z$ is $\mathbf{T}$-local iff $f^{*}:[B,Z] \to [A,Z]$ is a bijection for all $\mathbf{T}$-equivalences, $f: A \to B$, between connected cofibrant spaces, \\
	ii) if $Z = \sqcup_{i\in I} Z_i$, then $Z$ is $\mathbf{T}$-local iff $Z_i$ is $\mathbf{T}$-local for every $i$. In particular, a map of spaces which induces a bijection on connected components is a $\mathbf{T}$-localisation iff each component is a $\mathbf{T}$-localisation.
\end{Lemma}

\begin{proof}
	If $f: A \to B$ is a $\mathbf{T}$-equivalence between cofibrant spaces, then $f$ induces a bijection between the connected components of $A$ and $B$, so $f$ is a disjoint union of $\mathbf{T}$-equivalences $A_i \to B_i$, for $i$ in the set of connected components of $A$. Now, $[\sqcup A_i, Z] = \prod_i [A_i,Z]$, and $i)$ follows.  For $ii)$, if $A$ is connected we have $[A, \sqcup_i Z_i] = \sqcup_i [A, Z_i]$, and so $ii)$ follows from $i)$. 
\end{proof}

\begin{Lemma} \label{name55} 
	Let $f: X \to Y$ be a map of unbased spaces, with $X$ non-empty. Then the following are equivalent:
	
	i) $f$ is an unbased $\mathbf{T}$-localisation, \\
	ii) $f$ is a based $\mathbf{T}$-localisation for some $x \in X$, \\
	iii) $f$ is a based $\mathbf{T}$-localisation for all $x \in X$, \\
	iv) $f_+$ is a based $\mathbf{T}$-localisation, with respect to the adjoined basepoint $+$.
\end{Lemma}

\begin{proof}
	The key point is that if $Z$ is a $\mathbf{T}$-local based space, then it is also $\mathbf{T}$-local as an unbased space. This is a consequence of the fact that unbased homotopy classes $[A,Z]$ are equivalent to based homotopy classes $[A_+, Z]$, and the observation that if $A \to B$ is a $\mathbf{T}$-equivalence between cofibrant unbased spaces, then $A_+ \to B_+$ is a $\mathbf{T}$-equivalence between cofibrant based spaces. Now, $iii) \implies ii)$ is trivial, and $ii) \implies i)$ follows from the above. For $i) \implies iii)$, let $x \in X$. Since $\mathbf{T}$-localisations are preserved by composing with weak equivalences, we can assume that $X$ is a CW-complex and $f$ is a cofibration. Let $f_{\mathbf{T}}: X \to X_{\mathbf{T}}$ be a based $\mathbf{T}$-localisation. Then $f_{\mathbf{T}}$ is also an unbased $\mathbf{T}$-localisation, since $ii) \implies i)$. Therefore, there is a weak equivalence, $g$, such that $gf = f_{\mathbf{T}}$, and so $f$ is also a based $\mathbf{T}$-localisation, as desired. The fact that $iv) \implies i)$ follows from $ii) \implies i)$ and Lemma \ref{name54}ii), and $i) \implies iv)$ follows from Lemma \ref{name54}ii) and $i) \implies iii)$.  
\end{proof}

Returning to the equivariant setting, we have the following consequence:

\begin{Theorem} Let $\mathbf{T}$ be a localisation system. \\
	i) if $f: X \to Y$ is a based $\mathbf{T}$-localisation, then it is also an unbased $\mathbf{T}$-localisation, \\
	ii) if $f: X \to Y$ is a map of unbased $G$-spaces, then $f$ is a $\mathbf{T}$-localisation iff $f_+$ is a based $\mathbf{T}$-localisation. Moreover, if $X^G$ is non-empty, then $f$ is a $\mathbf{T}$-localisation iff $f$ is a based $\mathbf{T}$-localisation with respect to any $G$-fixed basepoint iff $f$ is a based $\mathbf{T}$-localisation with respect to all $G$-fixed basepoints.
\end{Theorem}

\begin{Remark}
	In this section, we have constructed functors which act as localisations at different homology theories at different fixed point spaces. Alternatively, as the referee noted, it might be possible to define a functor which acts as localisation at a set of primes at some fixed point spaces, and Bousfield-Kan's $R$-completion at others (for suitable $R$, \cite{BK72}). Of course, such a functor would no longer be idempotent but would instead lead to some very exotic $G$-spaces as `partial completions' of the original one.
\end{Remark}

\subsection{An algebraic analogue}

Before moving on to the theory of nilpotent $G$-spaces, we record the following result, which can be viewed as an algebraic analogue of the above theory. Recall that coefficient systems are functors $h \mathcal{O}^{op} \to \textbf{Ab}$, and there are free coefficient systems defined by:

\begin{Definition}
	The free coefficient system associated to the object $\orb{H}$ is defined by $\mathbf{F}_{\orb{H}}(\orb{K}) = \bigoplus_{h\mathcal{O}(\orb{K},\orb{H})} \mathbb{Z}$ along with the evident definition on morphisms. 
\end{Definition}

The free coefficient systems have the property that $Hom_{[h\mathcal{O}^{op},\mathbf{Ab}]}(A \otimes \mathbf{F}_{\orb{H}}, \mathbf{L}) \cong Hom_{\mathbf{Ab}}(A, \mathbf{L}(\orb{H}))$, where $A$ is any abelian group.

\begin{Lemma}
	Let $\mathbf{T}$ be a localisation system and let $\mathbf{A}$ and $\mathbf{B}$ be coefficient systems such that:
	
	i) if the coefficient of $\mathbf{T}(\orb{H})$ is 0, then $\mathbf{A}(\orb{H}) \otimes \mathbb{Z}_{T(\orb{H})} = 0$ and $\mathbf{B}(\orb{H})$ is $T(\orb{H})$-local, \\
	ii) if the coefficient of $\mathbf{T}(\orb{H})$ is 1, then $\mathbf{A}(\orb{H})$ is a $\mathbb{Z}[T(\orb{H})^{-1}]$-module and $\mathbf{B}(\orb{H})$ is $T(\orb{H})$-complete.
	
	Then $Ext^i_{[h\mathcal{O}^{op}, \mathbf{Ab}]}(\mathbf{A},\mathbf{B}) = 0$ for all $i \geq 0$.	
\end{Lemma}

\begin{proof}
	
	We first claim that if $\mathbf{T}(\orb{H})$ has coefficient 0, and $n$ is a product of primes not in $T(\orb{H})$, then $Ext^i_{[h\mathcal{O}^{op}, \mathbf{Ab}]}(\mathbf{F}_{\orb{H}} \otimes \sfrac{\mathbb{Z}}{n\mathbb{Z}},\mathbf{B}) = 0$ for all $i \geq 0$. The category $[h\mathcal{O}^{op}, \mathbf{Ab}]$ has enough injectives, \cite[Exercise 2.3.7]{W94}, so we can calculate this by taking an injective resolution $\{\mathbf{Q}_i\}$ of $\mathbf{B}$. Such a resolution is, in particular, an objectwise injective resolution of $\mathbf{B}(\orb{H})$, and $Hom_{[h\mathcal{O}^{op}, \mathbf{Ab}]}(\mathbf{F}_{\orb{H}} \otimes \sfrac{\mathbb{Z}}{n\mathbb{Z}}, \mathbf{Q}_i) = Hom_{\mathbf{Ab}}(\sfrac{\mathbb{Z}}{n \mathbb{Z}}, \mathbf{Q}_i (\orb{H}))$, so taking homology calculates $Ext^i(\sfrac{\mathbb{Z}}{n \mathbb{Z}}, \mathbf{B}(\orb{H}))$, which vanishes by the non-equivariant case. Similarly, if $\mathbf{T}(\orb{H})$ has coefficient 1, then $Ext^i_{[h\mathcal{O}^{op}, \mathbf{Ab}]}(\mathbf{F}_{\orb{H}} \otimes \mathbb{Z}[\mathbf{T}(\orb{H})^{-1}],\mathbf{B}) = 0$ by \cite[10.1.22]{MP12}.
	
	We will use this to define a $Hom_{[h\mathcal{O}^{op}, \mathbf{Ab}]}(-,\mathbf{B})$-acyclic resolution, $\{\mathbf{P}_i\}$, of $\mathbf{A}$. If $\mathbf{T}(\orb{H})$ has coefficient 0, there is a coproduct, $\mathbf{K}_{\orb{H}}$, of functors of the form $\mathbf{F}_{\orb{H}} \otimes \sfrac{\mathbb{Z}}{n\mathbb{Z}}$, with $n$ being a product of primes not in $T(\orb{H})$, such that there is a natural transformation $\mathbf{K}_{\orb{H}} \to \mathbf{A}$ which is a surjection at $\orb{H}$. If $\mathbf{T}(\orb{H})$ has coefficient 1, then there is a coproduct, $\mathbf{K}_{\orb{H}}$, of functors of the form $\mathbf{F}_{\orb{H}} \otimes \mathbb{Z}[\mathbf{T}(\orb{H})^{-1}]$, such that there is a natural transformation $\mathbf{K}_{\orb{H}} \to \mathbf{A}$ which is a surjection at $\orb{H}$. We define $\mathbf{P}_0 := \bigoplus_{\orb{H}} \mathbf{K}_{\orb{H}}$, so we have a surjection $\mathbf{P}_0 \to \mathbf{A}$, and $\mathbf{P}_0$ is $Hom_{[h\mathcal{O}^{op}, \mathbf{Ab}]}(-,\mathbf{B})$-acyclic by the previous paragraph.
	
	The key point now is that the functor $\mathbf{P}_0$ satisfies the conditions in i) and ii) that $\mathbf{A}$ does, and this follows from the fact that $\mathbf{T}$ is a localisation system. In more detail, $\mathbf{F}_{\orb{H}}(\orb{K})$ is only non-zero when there is a map $\orb{K} \to \orb{H}$ in $\mathcal{O}$, and then we have the following observations:
	
	i) if $S \subset T$, then a torsion group with no $T$-torsion is also a torsion group with no $S$-torsion, \\
	ii) if $S \subset T$, then a $\mathbb{Z}[T^{-1}]$-module is a $\mathbb{Z}[S^{-1}]$-module, \\
	iii) a torsion group with no $T$-torsion is a $\mathbb{Z}[T^{-1}]$-module.
	
	Therefore, we can inductively construct a $Hom_{[h\mathcal{O}^{op}, \mathbf{Ab}]}(-,\mathbf{B})$-acyclic resolution $\{\mathbf{P}_i\}$ of $\mathbf{A}$, since the kernel of $\mathbf{P}_0 \to \mathbf{A}$ also satisfies i) and ii) in the statement of the lemma. Using the first paragraph of the proof, we can use this acyclic resolution to compute $Ext^i_{[h\mathcal{O}^{op}, \mathbf{Ab}]}(\mathbf{A},\mathbf{B}) = 0$ for all $i \geq 0$, as desired. 	
\end{proof}

\section{Nilpotent $G$-spaces} \label{name187}

\subsection{The main theorems} \label{misc220}

We now move on to the theory of nilpotent $G$-spaces and we begin with the definition of a nilpotent $G$-space. This differs from the definition given in \cite[Ch. II]{M96} in that we do not require a common bound on the nilpotency classes at each fixed point space. To understand this, we will show in Subsection \ref{name80} that any nilpotent $G$-space can be approximated by a weak Postnikov tower, but if we assume a common bound on the nilpotency classes, then a nilpotent $G$-space can be approximated by a (strict) Postnikov tower, a distinction which becomes important when using co-HELP to deduce theorems about nilpotent spaces, as in \cite[Section 3.3]{MP12}. 

\begin{Definition} \label{name21}
	A based $G$-space $X$ is said to be nilpotent if $X^H$ is a nilpotent space for all subgroups $H$ of $G$.
\end{Definition}

In the unbased context, we have the following definition:

\begin{Definition}
	An unbased $G$-space $X$ is said to be componentwise nilpotent if for every subgroup $H$ of $G$, every component of $X^H$ is a nilpotent space.
\end{Definition}

In general, if we speak about componentwise nilpotent $G$-spaces we are working in an unbased context, and if we speak about nilpotent $G$-spaces we are working in a based context.

The non-equivariant theory of localisations and completions of nilpotent spaces was developed by Bousfield and Kan in \cite{BK72}. By reduction to fixed point spaces, we can immediately deduce one of the most important properties of localisations of componentwise nilpotent $G$-spaces:

\begin{Theorem}
	Let $\mathbf{T}$ be a localisation system where all the coefficients are 0. Let $f: X \to Y$ be a map from a componentwise nilpotent $G$-space $X$ to a $\mathbf{T}$-local unbased $G$-space $Y$, such that for every $H \leq G$, $f^H$ induces a bijection on connected components. Then, the following are equivalent:
	
	i) $f$ is a $\mathbf{T}$-localisation, \\
	ii) for all $H \leq G$, $* \geq 1$, and $b \in X^H$, $f^H_* : \pi_*(X^H,b) \to \pi_*(Y^H,f^H(b))$ is a $\mathbf{T}(\orb{H})$-localisation of nilpotent groups, \\
	iii)  for all $H \leq G$ and $* \geq 1$, $f^H_* : H_*(X^H) \to H_*(Y^H)$ is a direct sum of $\mathbf{T}(\orb{H})$-localisations, where the sum ranges over the connected components of $X^H$.
\end{Theorem}	

\begin{proof}
	This follows from \cite[Theorem 6.1.2]{MP12}, as well as Lemma \ref{name54}ii).
\end{proof}

Recall that if $T$ is a set of primes and $A$ is an abelian group, then $\mathbb{E}_T A$ and $\mathbb{H}_T A$ denote the zeroth and first derived functors of the classical $T$-completion of groups applied to $A$, respectively, \cite[Proposition 10.1.17]{MP12}. These functors can be extended to take nilpotent groups as input by using the homotopy groups of completions of Eilenberg-MacLane spaces. In the current context, we use the above definition of $\mathbb{E}_{\mathbf{T}} G$ and $\mathbb{H}_{\mathbf{T}} G$ for sets of primes with coefficient 1. If, instead, $\mathbf{T}$ is a set of primes with coefficient 0, and $G$ is a nilpotent group, we define $\mathbb{E}_{\mathbf{T}} G = G_T$ and $\mathbb{H}_{\mathbf{T}} G = 0$. This corresponds to using the homotopy groups of localisations of Eilenberg-MacLane spaces. A system of nilpotent groups, $\underline{G}$, is a continuous functor from $\mathcal{O}^{op}$ to the category of nilpotent groups, and we call such a system $\mathbf{T}$-local if it is pointwise $\mathbf{T}(\orb{H})$-local. The $\mathbf{T}$-localisation $K(\underline{G},1) \to K(\underline{G},1)_{\mathbf{T}}$ specifies a homomorphism $\underline{G} \to \mathbb{E}_{\mathbf{T}}(\underline{G})$ and the, up to homotopy, universal property of $\mathbf{T}$-localisation implies the following universal property:  

\begin{Lemma}
	Let $\underline{G}$ and $\underline{H}$ be systems of nilpotent groups, with $\underline{H}$ $\mathbf{T}$-local. Then any homomorphism $f:\underline{G} \to \underline{H}$ factors uniquely through the $\mathbf{T}$-localisation $\underline{G} \to \mathbb{E}_{\mathbf{T}} \underline{G}$. 
\end{Lemma}

\begin{proof}
	This follows from the fact that $[K(\underline{G},1)_{\mathbf{T}}, K(\underline{H},1)] \cong [K(\underline{G},1),K(\underline{H},1)]$. 
\end{proof}

If $X$ is a nilpotent $G$-space, then this universal property defines a map from $\mathbb{E}_{\mathbf{T}} \underline{\pi}_i (X) \to \underline{\pi}_i(X_{\mathbf{T}})$ and we have the following theorem:

\begin{Theorem} \label{name12} 
	If $X$ is a nilpotent $G$-space, then there is a natural short exact sequence:
	
	\begin{center}
		$1 \to \mathbb{E}_{\mathbf{T}} \underline{\pi}_i (X) \to \underline{\pi}_i (X_{\mathbf{T}}) \to \mathbb{H}_{\mathbf{T}} \underline{\pi}_{i-1} (X)  \to 1$
		
	\end{center}
	
	If $f: X \to Y$ is a map between componentwise nilpotent $G$-spaces such that each $f^H$ induces a bijection on connected components, and $\mathbb{H}_{\mathbf{T}(\orb{H})}(\pi_i(X^H,x)) = 0$ for all $H \leq G, i \geq 1$ and $x \in X^H$, then the following are equivalent:
	
	i) $f$ is a $\mathbf{T}$-localisation, \\
	ii) for all $i \geq 1$, $H \leq G$ and $x \in X^H$, $\pi_i(X^H,x) \to \pi_i(Y^H,f^H (x))$ is a $\mathbf{T}(\orb{H})$-localisation.
	
	For example, the hypothesis holds if, for all $H$, $X_{\mathbf{T}(\orb{H})}^H$ is a disjoint union of f$\mathbb{Z}_{T(\orb{H})}$-nilpotent spaces.
	
\end{Theorem}

\begin{proof}
	This follows from \cite[Theorem 11.1.2, Proposition 10.1.23]{MP12}, as well as Lemma \ref{name54}.
\end{proof}

Non-equivariantly, the fact that $Ext(\mathbb{H}_T B, \mathbb{E}_T A) = 0$, \cite[Corollary 10.4.9]{MP12}, implies that the short exact sequence of Theorem \ref{name12} splits, however, equivariantly the sequence does not necessarily split as the following example shows:

\begin{Example} Take $G = C_2$. Then, consideration of Elmendorf's theorem, \cite[Theorem 1]{E83}, shows that to find a counterexample to the splitting, we can use the following counterexample to the naturality of the splitting in the non-equivariant case.  For this, we let $X =  K(\frac{\mathbb{Z}[p^{-1}]}{\mathbb{Z}},1)$, so that $\widehat{X}_p =  K(\widehat{\mathbb{Z}}_p,2)$, and a map $X \to K(\widehat{\mathbb{Z}}_p,2)$ is equivalent to a homomorphism $\widehat{\mathbb{Z}}_p  \to \widehat{\mathbb{Z}}_p$. Then, any non-zero homomorphism, such as the identity, suffices to show that the splitting cannot be natural. \end{Example}

We close with a few examples of nilpotent $G$-spaces:

\begin{Example}
	If $S^V$ is a representation sphere, then all fixed point spaces are spheres, so $S^V$ is componentwise nilpotent. In fact, $S^V$ is nilpotent iff $V$ has a non-zero fixed point.
\end{Example}

\begin{Example}
	We know, from \cite[Theorem 10]{LM86}, that if $\Pi$ is a normal subgroup of a compact Lie group $\Gamma$, $G = \Gamma / \Pi$ and $H$ is a subgroup of $G$, then $B(\Pi; \Gamma)^H = \sqcup B(\Pi \cap N_{\Gamma} \Lambda)$, where the disjoint union runs over the $\Pi$-conjugacy classes of subgroups $\Lambda$ of $\Gamma$ such that $\Lambda \cap \Pi = 1$ and $q(\Lambda) = H$, where $q$ is the quotient map. Therefore, if $\Pi$ is abelian, then $B(\Pi; \Gamma)$ is a componentwise nilpotent $G$-space. In this case, the effect of $\mathbf{T}$-localisation on fixed point spaces will be a disjoint union of finite products of localisations/completions of $K(\mathbb{Z},2)$ and $K(A,1)$, where $A$ is an abelian group.
\end{Example}

\begin{Remark}	
We are working in this section with nilpotent $G$-spaces, but it is clear that our general approach also allows us to port non-equivariant results on localisations of virtually nilpotent spaces into the equivariant context. For example, the real projective spaces $\mathbb{RP}^n$ are virtually nilpotent spaces, \cite[2.6]{DFK77}, and so, as for all virtually nilpotent spaces, we can fracture them into complete and rational parts, \cite[Theorem 4.1]{DFK77}. We will discuss such fracture theorems more in Subsections \ref{SS32} and \ref{name182}
\end{Remark}

\subsection{$T$-localisation and fibre squares} \label{misc221}

From now on, we work in a based context. In this subsection, we discuss how $\mathbf{T}$-localisation interacts with fibre sequences and homotopy pullbacks. Let $N(f,g)$ denote the double mapping path space of $f$ and $g$. We have:

\begin{Theorem} \label{name19}
	Let $f:X \to Z$ and $g:Y \to Z$ be maps of nilpotent $G$-spaces such that $N(f,g)$ is $G$-connected. If we have a commutative diagram:
	
	\[	\begin{tikzcd}
		X \arrow{r}{f} \arrow{d} & Z \arrow{d} & Y \arrow[swap]{l}{g} \arrow{d} \\
		X^{'} \arrow[swap]{r}{f^{'}} & Z^{'} & Y^{'} \arrow{l}{g^{'}} \\
	\end{tikzcd}
	\]
	
	such that the vertical maps are $\mathbf{T}$-localisations, then the induced map $N(f,g) \to N(f^{'},g^{'})$ is a $\mathbf{T}$-localisation.
\end{Theorem}

\begin{proof}
	This follows from \cite[Theorem 3.2]{R24}.
\end{proof}

It follows that any functorial $\mathbf{T}$-localisation, such as one obtained via the small object argument applied to the model structure of Theorem \ref{name8}, preserves homotopy fibre squares of nilpotent $G$-spaces. The special case where $Y = Y^{'} = *$ results in the connected fibre lemma.

\subsection{Fracture theorems for virtually nilpotent $G$-spaces} \label{SS32} 

In this subsection, we move on to derive fracture squares associated to nilpotent $G$-spaces. In the non-equivariant setting, we have fracture squares relating to localisation and completion, first studied by Bousfield and Kan in \cite[V.6, VI.8]{BK72} but see also \cite[Theorem 8.1.3, Theorem 13.1.4]{MP12}, and we would like to generalise these results to the equivariant setting, perhaps localising and completing at different sets of primes at each fixed point space. For example, the following two squares are homotopy fibre squares associated to a $T$-local nilpotent space $X$, where $T$ is a set of primes containing 7:

\begin{equation} \label{dia1}
	\begin{tikzcd}
		X \arrow{r} \arrow{d} & \widehat{X}_T \arrow{d} \\
		X_0 \arrow{r} & (\widehat{X}_T)_0 \\
	\end{tikzcd}
\end{equation}

\begin{equation} \label{dia2}
	\begin{tikzcd}
		X \arrow{r} \arrow{d} & (\prod_{p \in T \backslash \{7\}} X_p) \times \widehat{X}_7 \arrow{d} \\
		X_0 \arrow{r} & (\prod_{p \in T \backslash \{7\}} X_p)_0 \times (\widehat{X}_7)_0 \\
	\end{tikzcd}
\end{equation}

In \ref{dia2}, we complete at 7 to illustrate to point that there are an abundance of fracture squares that we can ask for, especially in the equivariant case. With this in mind, the following theorem subsumes all of the examples that we are aware of. We focus on nilpotent $G$-spaces for now but the proof works for virtually nilpotent $G$-spaces ($G$-spaces whose fixed point spaces are virtually nilpotent), as we will see in the next subsection.

\begin{Theorem} \label{name88}
	Consider a commutative square of virtually nilpotent $G$-spaces:
	\begin{center}

		\begin{tikzcd}
			X \arrow{r}{f} \arrow[swap]{d}{\phi} & Y \arrow{d}{\psi} \\
			A \arrow[swap]{r}{g} & B \\
		\end{tikzcd}
	\end{center}	
	
	such that, for each subgroup $H$ of $G$, there are sets of primes $T_H, S_H$ satisfying:
	
	i) $X^H,Y^H$ are $T_H$-local and $A^H,B^H$ are $S_H$-local, \\
	ii) $f^H$ is an $\mathbb{F}_{T_H}$-equivalence and $g^H$ is an $\mathbb{F}_{S_H}$-equivalence, \\
	iii) $\phi^{H}, \psi^{H}$ are $\mathbb{Q}$-equivalences.
	
	Then the square is a homotopy fibre square.

\end{Theorem}

\begin{proof}
	Since taking fixed points detects homotopy fibre squares, we can reduce the theorem to a pointwise statement with sets of primes $T$ and $S$. In this case, the theorem follows from (\ref{dia1}) and successive applications of the pasting lemma for homotopy pullbacks, \cite[Proposition 13.3.15]{H03}.
\end{proof}

\subsection{An introduction to virtually nilpotent spaces} \label{name182}

In this subsection, we introduce the theory of virtually nilpotent spaces. The theory was originally developed by Dror Farjoun, Dwyer and Kan in \cite{DFK77}, and the treatment here follows along similar lines, the main difference being that we work solely with homological localisations.  Despite working in a non-equivariant context, we will continue to use notation consistent with the previous sections. In particular, if $T$ is a set of primes, we let $X_T$ denote the $H(-;\mathbb{Z}_T)$-localisation of $X$ and $X_{\mathbf{T}}$ denote the $H(-; \oplus_{p \in T} \mathbb{F}_p)$-localisation of $X$. We also remind the reader that all (virtually) nilpotent spaces are assumed to be connected.

We begin with the definition of a $\mathbb{Z}_T$-perfect group, but first recall that a $T$-number is a natural number which is a product of primes not in $T$:

\begin{Definition}
	If $G$ is a group and $T$ is a set of primes, define $\Gamma_T^i G$ to be the subgroup of $G$ consisting of those $g$ for which there exists a $T$-number $r$ such that $g^r \in \Gamma^iG$. A group is said to be $\mathbb{Z}_T$-perfect if $\Gamma_T^1 G = G$ or, equivalently, $H_1(G; \mathbb{Z}_T) = 0$. We call a space, $Y$, $\mathbb{Z}_T$-perfect if $H_1(Y; \mathbb{Z}_T) = 0$.
\end{Definition}

Note that $[G, \Gamma^i_T G] \subset \Gamma^{i+1}_T G$ since if $g \in G$, then in the subgroup, $A$, generated by $g$ and $\Gamma^i_T G$, $\Gamma^j A$ ($j \geq 1$) is generated by $j$-fold commutators containing at least one element of $\Gamma^i_T G$, by \cite[Lemma 2.6]{CMZ17}. So, inductively, $\Gamma^{i+1} A, ..., \Gamma^{1} A \subset \Gamma_T^{i+1} G$, using the epimorphisms of \cite[Corollary 2.10]{CMZ17}. So $\{\Gamma^i_T G\}$ is a bona fide central series which, of course, coincides with the lower central series in the case where $G$ is a $T$-local nilpotent group.

\begin{Lemma} \label{name174}
	If $1 \to K \to G \to H \to 1$ is an exact sequence of groups with $K$ a nilpotent group of class $c$ which is $p$-divisible for $p \notin T$, and $H$ is $\mathbb{Z}_T$-perfect, then $\Gamma_T^c G$ is $\mathbb{Z}_T$-perfect.
\end{Lemma}

\begin{proof}
	We first show that, for arbitrary $i$, any $g \in G$ can be written in the form $kx$ with $k \in K$ and $x \in \Gamma_T^i G$. To see this note that $\Gamma^i_T H = H$ since $H$ is $\mathbb{Z}_T$-perfect and $[H, \Gamma^j_T H] \subset \Gamma^{j+1}_T H$. So for any $g \in G$ there exists a $T$-number, $r$, with $g^r = ay$ with $a \in K$ and $y \in \Gamma^i G$. Inductively, we may assume $g \in \Gamma^{i-1}_T G$. Since $K$ is $p$-divisible for $p \notin T$, we can write $a^{-1} = b^r$ for $b \in K$. Then, using the fact that $[G, \Gamma^{i-1}_T G] \subset \Gamma^{i}_T G$, we see that $(bg)^r = b^r g^r x = yx$ where $x \in \Gamma^i_T G$ and so $bg \in \Gamma^i_T G$, which yields the required result.
	
	Therefore, by \cite[Lemma 2.6]{CMZ17}, $\Gamma^c G$ is generated by commutators $[g_0,...,g_c]$ where each $g_i$ is an element of either $K$ or $\Gamma^{c+1}_T G$ and at most one of the $g_i$ is not in $K$, since $K$ has nilpotency class $c$. It follows that $\Gamma^{c}_T G = \Gamma^{c+1}_T G$. Moreover, $\Gamma^{c+1} G$ is generated by commutators with at least two terms not in $K$ and so $\Gamma^c_T G = \Gamma^{c+1}_T G \subset \Gamma^1_T (\Gamma_T^c G)$. Therefore, $\Gamma^c_T G$ is $\mathbb{Z}_T$-perfect.
\end{proof}

Let $n \geq 1$. Recall that a map $f: X \to Y$ between connected spaces is called an $n$-equivalence if $\pi_i(Ff) = 1$ for $0 \leq i \leq n-1$ and a homology $n$-equivalence if $\tilde{H}_i (Cf) = 0$ for $0 \leq i \leq n$. For the purposes of the current discussion, we keep track of the connectivity of the maps involved in the upcoming lemmas - this will eventually be used to conclude that certain spaces are still virtually nilpotent after passing to inverse limits. However, we note that in \cite{DFK77}, the relative connectivity lemma of \cite[IV.5.1]{BK72} was used to derive such results. We begin with \cite[5.2]{DFK77}:

\begin{Lemma} \label{name173}
	For $i = 1,2$, let $N_i \to X_i \to Y_i$ be a fibre sequence of connected spaces such that $N_i$ is nilpotent and $Y_i$ is $\mathbb{Z}_T$-perfect. Then each $(X_i)_T$ is nilpotent. Suppose that we have a commutative square:	
	\[\begin{tikzcd}
		X_1 \arrow{r}{f} \arrow{d} & X_2 \arrow{d} \\
		Y_1 \arrow[swap]{r}{g} & Y_2 
	\end{tikzcd}
	\]
	
	If $Y_1 = Y_2 := Y$, $g = 1$ and $f$ is an $n$-equivalence, then $f_T$ is an $n$-equivalence. If, instead, the induced map $N_1 \to N_2$ between homotopy fibres is a $1$-equivalence, then $f_T$ is a $1$-equivalence.
\end{Lemma}

\begin{proof}
	Localizing the fibre, \cite[Theorem 6.3.1]{H03}, we obtain fibre sequences $(N_i)_T \to X_i^{'} \to Y_i$. The Serre spectral sequence implies that the maps $X_i \to X_i^{'}$ are $\mathbb{Z}_T$-equivalences. In the $g = 1$ case, since each $N_i$ is nilpotent, $X_1^{'} \to X_2^{'}$ is still an $n$-equivalence. Similarly, if $N_1 \to N_2$ is a $1$-equivalence, so is $(N_1)_T \to (N_2)_T$.  
	
	Let $K_i = \pi_1((N_i)_T)$, $G_i = \pi_1(X_i^{'})$ and $H_i = \pi_1(Y_i)$ so that we have an exact sequence $K_i \to G_i \to H_i \to 1$ with $H_i$ $\mathbb{Z}_T$-perfect and $K_i$ nilpotent and $T$-local. Lemma \ref{name174} implies that $\Gamma^c_T G_i$ is $\mathbb{Z}_T$-perfect, where $c$ is the maximum of the nilpotency classes of $K_1$ and $K_2$, so consider the diagrams of fibrations:
	
	\begin{equation} \label{name175}
		\begin{tikzcd} 
			P_i \arrow{r} \arrow{d} & (N_i)_T \arrow{d} \arrow{r} & K(G_i / \Gamma_T^c G_i, 1) \arrow{d} \\
			Q_i \arrow{r} \arrow{d} & X_i^{'} \arrow{d} \arrow{r} & K(G_i / \Gamma_T^c G_i, 1) \arrow{d} \\
			Y_i \arrow{r} & Y_i \arrow{r} & * \\
		\end{tikzcd} 
	\end{equation}

	We saw in the proof of Lemma \ref{name174} that $K_i \to G_i / \Gamma^c_T G_i$ is surjective, so all the spaces are connected. Note that in a situation such as \ref{name175}, $\pi_1(X_i{'})$ acts, up to homotopy, on all of the spaces in the diagram, including $P_i$ since we can identify $P_i$ with the homotopy fibre of the map from $X_i^{'}$ to the homotopy pullback of the maps $K(G_i / \Gamma^c_T G_i, 1) \to *$ and $Y_i \to *$. (We do not require the target of these maps to be $*$, but it is in this instance.) The induced action of $K_i$ on $P_i$ can be identified with the action defined by top horizontal fibre sequence. Moreover, as in \cite[Lemma 4.2.7]{R23}, $\pi_1(X_i^{'})$ acts on the Serre spectral sequence of the fibration $P_i \to Q_i \to Y_i$. Using these observations, as in \cite[Lemma 5.1]{DFK77}, a routine Serre spectral sequence argument shows that $G_i/ \Gamma_T^c G_i$ acts nilpotently on $H_*(Q_i; \mathbb{Z}_T)$.
	
	Since $\Gamma^c_T G_i$ is $\mathbb{Z}_T$-perfect, $H_1(Q_i; \mathbb{Z}_T) = 0$ and so $(Q_i)_T$ is simply connected by a plus-construction style argument, \cite[Proposition 4.40]{H02}. Therefore, by localising fibres, we obtain quasi-nilpotent fibrations $(Q_i)_T \to X_i^{''} \to K(G_i/\Gamma^c_T G_i,1)$ with nilpotent fibre. This implies that the fibrations are nilpotent via \cite[Corollary 2.2]{H76}. It follows that $X_i^{''}$ is nilpotent which implies that $(X_i^{''})_T \simeq (X_i)_T$ is nilpotent, as desired.
	
	 In the $g = 1$ and $n \geq 2$ case, $Q_1 \to Q_2$ is an $n$-equivalence, so also a homology $n$-equivalence, and so $(Q_1)_T \to (Q_2)_T$ is a homology $n$-equivalence between simply connected spaces, so is also an $n$-equivalence. It follows that $X^{''}_1 \to X^{''}_2$ is an $n$-equivalence, and, therefore, so is $f_T$. For the remaining cases, note that $\pi_1((X_i)_T) \cong (G_i / \Gamma^c_T G_i)_T$ and, since $K_i$ surjects onto $G_i / \Gamma^c_T G_i$ and $K_1 \to K_2$ is surjective, it follows that $f_T$ is a $1$-equivalence.
\end{proof}

Next, we record the following straightforward consequence of \cite[Proposition 12.7]{B75}, the Zeeman comparison theorem of \cite{HR76} and the observation that a finite $p$-group always acts nilpotently on an $\mathbb{F}_p$-module, \cite[pg. 215]{BK72}:

\begin{Lemma} \label{name176}
	If $F \to E \to B$ is a fibre sequence of connected spaces such that $\pi_1(B)$ is a finite $p$-group, then:
	
	i) $\pi_1(B) \to \pi_1(B_{\mathbf{p}})$ is an isomorphism, \\
	ii) $F_{\mathbf{p}} \to E_\mathbf{p} \to B_\mathbf{p}$ is a fibre sequence, \\
	iii) if $B$ and $F$ are $\mathbf{p}$-local, then $E$ is $\mathbf{p}$-local.
\end{Lemma}

In fact, the fundamental group of $X_{\mathbf{p}}$ (or $X_T$) is always determined by $\pi_1(X)$, by \cite[Lemma 7.3]{B75}. The next lemma helps us understand $\mathbf{p}$-localisations of nilpotent-by-finite spaces:

\begin{Lemma} \label{name178} For $i = 1,2$, let $N_i \to X_i \to Y_i$ be fibre sequences with $N_i$ nilpotent and $\pi_1 (Y_i)$ finite. Then, for all primes $p$, $(X_i)_\mathbf{p}$ is nilpotent-by-finite. Suppose that we have a commutative square:	
	\[\begin{tikzcd}
		X_1 \arrow{r}{f} \arrow{d} & X_2 \arrow{d} \\
		Y_1 \arrow[swap]{r}{g} & Y_2 
	\end{tikzcd}
	\]
	If $g = 1_Y$ and $f$ is an $n$-equivalence $(n \geq 1)$, then so is $f_\mathbf{p}$. If $g$ and the induced map $N_1 \to N_2$ are $1$-equivalences, so is $f_{\mathbf{p}}$.
\end{Lemma}

\begin{proof}
	Consider the diagrams:	
	\begin{equation} \label{name181}
		\begin{tikzcd}
			N_i \arrow{d} \arrow{r} & N_i \arrow{d} \arrow{r} & * \arrow{d} \\
			E_i \arrow{d} \arrow{r} & X_i \arrow{d} \arrow{r} & K(P_i,1) \arrow{d} \\
			F_i \arrow{r} & Y_i \arrow{r} & K(P_i,1) 
	\end{tikzcd} \end{equation}

	where $P_i$ is the quotient of $\pi_1(Y_i)$ by its maximal $\mathbb{Z}_p$-perfect subgroup, so $P_i$ is a finite $p$-group. Now, $(E_i)_\mathbf{p}$ is nilpotent by Lemma \ref{name173}, and $(E_i)_\mathbf{p} \to (X_i)_\mathbf{p} \to K(P_i,1)$ is a fibre sequence, by Lemma \ref{name176}. In the $g=1$ case, $E_1 \to E_2$ is an $n$-equivalence and $F_1 = F_2$. So $f_\mathbf{p}$ is an $n$-equivalence, using Lemma \ref{name173}. In the remaining case, $g$ is a $1$-equivalence, so $P_1 \to P_2$ is surjective. So the fact that $f_\mathbf{p}$ is a $1$-equivalence follows from Lemma \ref{name173} as well.
\end{proof}

	Note that if $X$ is $T$-local and $p \notin T$, then all homotopy groups of $X$ are uniquely $p$-divisible. This can be seen by considering the $\mathbb{Z}_T$-equivalences $S^n \to \vee S^n \to S^n$ which represent multiplication by $p$ on homotopy groups. Therefore, the case where $X$ is $T$-local and $p \notin T$ is an important special case of the following theorem:

\begin{Lemma} \label{name179}
	Let $N \to X \to K(F,1)$ be a fibre sequence of connected spaces with $N$ nilpotent and $F$ a finite group. If the homotopy groups of $X$ are uniquely $p$-divisible, then $X_\mathbf{p} \simeq *$.
\end{Lemma}

\begin{proof}
Since $F$ must be $p$-divisible, $p$ does not divide the order of $F$. Therefore, $F$ and the homotopy groups of $N$ are uniquely $p$-divisible. Since $N$ is nilpotent, this implies $N_{\mathbf{p}} = *$. Also, consideration of the transfer map shows that $\tilde{H}_*(F; \mathbb{F}_p) = 0$, \cite[pg.129, Problem 3]{M99}. Use of the Serre spectral sequence now implies that $\tilde{H}_*(X; \mathbb{F}_p) = 0$, as desired.
\end{proof}

We now deduce most of the following properties of virtually nilpotent spaces from the nilpotent-by-finite case. We say most because some of the results about $X_T$ follow from the arithmetic square theorem which we will derive shortly:

\begin{Theorem} \label{name180}
	If $X$ is virtually nilpotent, then so are $X_T, X_\mathbf{T}$ and $X_\mathbf{T} \simeq \prod_{p \in T} X_\mathbf{p}$. Also, $X_0$ is nilpotent. If $X$ is $T$-local and $p \notin T$, then $X_\mathbf{p} = *$. Moreover, if $f:X_1 \to X_2$ is an $n$-equivalence ($n \geq 1$) between virtually nilpotent spaces, then $f_\mathbf{T}, f_0$ are $n$-equivalences and $f_T$ is an $(n-1)$-equivalence if $n \geq 2$.
\end{Theorem}

\begin{proof}
	
	We begin with the connectivity results for $f$, assuming that $X_1$ and $X_2$ are nilpotent-by-finite. First note that $(X_i)_{\mathbf{T}} \simeq \prod (X_i)_{\mathbf{p}}$ by Lemma \ref{name179} and the Kunneth theorem. If $f$ is an $n$-equivalence and $n \geq 2$, then it is straightforward to construct a finite group $F$ and a commutative triangle as below such that the respective homotopy fibres are nilpotent:
	
	\[ \begin{tikzcd}
		X_1 \arrow{dr} \arrow{rr} & & X_2 \arrow{dl} \\
		& K(F,1) &
	\end{tikzcd}\]

	Therefore, $f_{\mathbf{T}}$ is an $n$-equivalence by Lemma \ref{name178} and $f_0$ is an $n$-equivalence by Lemma \ref{name173}. If $n = 1$, it is straightforward to construct a diagram:
	
	\[	 \begin{tikzcd}
		X_1 \arrow{d} \arrow{r} & X_2 \arrow{d} \\
		Y_1 := K(B_1,1) \arrow{r} & Y_2 := K(B_2,1) 
	\end{tikzcd} \]
	
	where $B_1$ and $B_2$ are finite groups, the homotopy fibres $N_1$ and $N_2$ are nilpotent and the map $N_1 \to N_2$ is a $1$-equivalence. Therefore, $f_{\mathbf{T}}$ is a $1$-equivalence by Lemma \ref{name178} and $f_0$ is a $1$-equivalence by Lemma \ref{name173}.

	If $X$ is virtually nilpotent, then using the connectivity results above, \cite[Proposition 2.5.9]{MP12} and \cite[Proposition 12.8]{B75}, we have $X \simeq lim P_n X \to  lim (P_n X)_\mathbf{T} \simeq X_{\mathbf{T}}$.  It follows from Lemma \ref{name178} that $(P_nX)_\mathbf{T}$ is virtually nilpotent since we also know that there is some $T^{'} \subset T$ containing all but finitely many primes of $T$ such that $(P_n X)_{T^{'}}$ is nilpotent, by Lemma \ref{name173}. Since $P_{n+1}X \to P_n X$ is an $(n+1)$-equivalence,  and using the algebraic characterisation of virtually nilpotent spaces \cite[2.2]{DFK77}, it is now clear that $X_\mathbf{T}$ is virtually nilpotent.  Similarly, use of the algebraic characterisation of nilpotent spaces and Lemma \ref{name173} shows that $X_0$ is nilpotent. If $X$ is $T$-local, then each $P_n X$ has uniquely $p$-divisible homotopy groups, so $X_\mathbf{p} = *$ by Lemma \ref{name179}. In fact, each $P_n X$ is also $T$-local by \cite[Theorem 5.5]{B75}.
	
	The fact that $X_T$ is virtually nilpotent and the connectivity results for $f_T$ follow from the arithmetic square theorem, using \cite[Corollary 2.2.3]{MP12}.
\end{proof}

\begin{Lemma} \label{name177}
	Let $n \geq 1$ and $F \to E \to K(Q,n)$ be a fibre sequence of connected spaces with $E$ virtually nilpotent and $Q$ a rational nilpotent group. Then:
	
	i) for all primes $p$, $F \to E$ is an $\mathbb{F}_p$-equivalence, \\
	ii) $F_0 \to E_0 \to K(Q,n)$ is a fibre sequence.
\end{Lemma}

\begin{proof}
	The universal property of Postnikov approximations, \cite[1.5]{H03}, implies that if $m \geq n$, then $P_m F \to P_m E \to K(Q,n)$ is a fibre sequence. So we can reduce to the case where $E$ is nilpotent-by-finite and consider a diagram:
	
	\[	\begin{tikzcd}
		R \arrow{r} \arrow{d} & N \arrow{r} \arrow{d} & K(Q,n) \arrow{d} \\
		F \arrow{r} \arrow{d} & E \arrow{r} \arrow{d} & K(Q,n) \arrow{d} \\
		K(B,1) \arrow{r} & K(B,1) \arrow{r} & * \\
	\end{tikzcd} \]
	
	where $N$ is nilpotent and $B$ is a finite group. If $n = 1$, then $Q / im(\pi_1 (N))$ is a divisible finite group, so is trivial. Therefore, $R$ is connected and $\mathbf{p}$-localising the top row, which is a nilpotent fibration, shows that $R \to N$ is an $\mathbb{F}_p$-equivalence. Therefore, the Serre spectral sequence implies that $F \to E$ is an $\mathbb{F}_p$-equivalence. Moreover, as in the proof of Lemma \ref{name173}, $Q$ acts nilpotently on $H_*(F;A)$, for all abelian groups $A$. Since $F_0$ is nilpotent, it follows, from \cite[Corollary 2.2]{H76} and localising the fibre, that $(-)_0$ preserves the fibre sequence $F \to E \to K(Q,n)$.
\end{proof}

Using Lemma \ref{name177} and the relative Postnikov towers described, for instance, in \cite[Theorem 2.3.11]{R23}, it is now straightforward to derive the second arithmetic square theorem of \cite[Theorem 4.4]{DFK77}. To deduce the first arithmetic square theorem from the second we need the following connectivity result, which is \cite[Proposition 3.5]{DFK77}:

\begin{Lemma} \label{name189}
	If $X$ is a virtually nilpotent space, then the function $\pi_1(X_0) \times \pi_1 (X_{\mathbf{T}}) \to \pi_1((X_{\mathbf{T}})_0)$, defined via composition with group multiplication, is surjective.
\end{Lemma}

\begin{proof}
	We may as well assume $X$ is nilpotent-by-finite and we have a fibre sequence $N \to X \to K(F,1)$ with $F$ finite and $N$ nilpotent. We can prove the nilpotent case of the lemma as in \cite[Lemma 3.3]{R24}, so it suffices to prove that $\pi_1((N_{\mathbf{T}})_0) \to \pi_1((X_{\mathbf{T}})_0)$ is surjective. Recall from Lemma \ref{name173} that if $K(F,1)$ is $\mathbb{Z}_T$-perfect, then $X_T$ is nilpotent, $\pi_1(X_T) = (G / \Gamma_T^c G)_T$ and $(\pi_1(N))_T \to G / \Gamma_T^c G$ is surjective. Similarly, given that $X_{\mathbf{T}} = \prod X_{\mathbf{p}}$, consideration of Diagram \ref{name181} shows that the image of $\pi_1(N_{\mathbf{T}})$ in $\pi_1(X_{\mathbf{T}})$ has finite index. We know that $X_{\mathbf{T}}$ is nilpotent-by-finite, which gives rise to another finite index nilpotent subgroup, $K$, of $\pi_1(X_{\mathbf{T}})$. Since $im(\pi_1(N_{\mathbf{T}})) \cap K$ has finite index in $K$, and $K_0$ surjects onto $\pi_1((X_{\mathbf{T}})_0)$, we obtain the desired result.
\end{proof}

The first arithmetic square theorem, which is the analogue of \ref{dia1}, can now be derived as on \cite[pg. 247]{DFK77}. It is now clear that the proof of Theorem \ref{name88} also holds in the virtually nilpotent context. So, for example, if $X$ is a virtually nilpotent space, we also have a `local' fracture square:

\[
\begin{tikzcd}
	X_T \arrow{r} \arrow{d} & \prod_{p \in T} X_p \arrow{d} \\
	X_0 \arrow {r} & (\prod_{p \in T} X_p)_0 	
\end{tikzcd}
\]

\subsection{Fracture theorems for homotopy classes} \label{name196}

When the underlying set of primes of $\mathbf{T}$ is constant, we can derive a fracture square for homotopy classes $[K,X]$, under certain finiteness hypotheses on $K$ and $X$. In particular, $K$ will always be a finite \textit{based} $G$-CW complex, by which we mean that $K$ is built by starting with a $G$-fixed basepoint and inductively attaching finitely many cells along based maps out of $G$-spaces of the form $(\frac{G}{H})_+ \wedge S^n$ with image in the $n$-skeleton. In order to give what we feel is the cleanest exposition of our main result, and the corresponding counterexample when the underlying set of primes of $\mathbf{T}$ is not constant, we begin by recalling some preliminaries on homotopy pullbacks.

\begin{Definition}
	Let $f: K \to X$ be a map of $G$-spaces. We define $[K \wedge (I^n)_+, f]$ (technically, \\ $[(K) \wedge ((I^{ \times n})_+), f]$) to be the set of homotopy classes of maps $K \wedge (I^n)_+ \to X$ relative to the boundary $K \wedge (\partial I^n)_+$, where at each point on the boundary $\partial I^n$, the induced map is equal to $f$.
\end{Definition}

We will use the following lemma, which is closely related to the results of \cite[V.5]{BK72}:

\begin{Lemma} \label{name95} Let $A \xrightarrow{i} B \xrightarrow{j} Ci$ be a cofibre sequence of $G$-spaces, and $f: Ci \to X$ a map of $G$-spaces. Then there is a natural long exact sequence of groups:
	
	\[ 
	... \to [B \wedge (I^2)_+, fj] \to [\Sigma^2 A, X] \to [Ci \wedge I_+, f] \to [B \wedge I_+, fj] \to [\Sigma A, X]
	\]
	
	Moreover, the image of $[\Sigma^2 A, X]$ in $[Ci \wedge I_+,f]$ is central.
	
\end{Lemma}

\begin{proof}
	Modify $f$ so that it is radially constant in a neighbourhood of the boundary of the cone. Consider the sequence of based maps:
	
	\[
	... \to \Omega Map_G(A,X) \xrightarrow{\partial} Map_G(Ci,X) \xrightarrow{j^*} Map_G(B,X) \xrightarrow{i^*} Map_G(A,X) 
	\]
	
	where the spaces are given basepoints $f$ and the constant loop to $f$. The fact that $f$ is radially constant in a neighbourhood of the boundary of the cone allows us to define a based map, which is also a weak equivalence, $Fi^* \to Map_G(Ci,X)$. The map $\partial$ is then induced by a comparison of the fibre sequences associated to $j^*$ and $Fi^* \to Map_G(B,X)$. It follows that $[S^1,-]$ takes the above sequence of maps to an exact sequence of groups, since it does so for the homotopy fibre sequence induced by $i^*$. The fact that the image of $\pi_1(\partial)$ is central follows from \cite[Lemma 1.4.7 v).]{MP12}.	 
\end{proof}

Let $N(f,g)$ denote the double mapping path space associated to maps $f:X \to A$ and $g: Y \to A$. We will make use of the following result on homotopy classes of maps into a homotopy pullback:

\begin{Lemma} \label{name89}
	
	Let $K$ be a based $G$-CW complex. Then the natural map of pointed sets:	
	\[
	[K,N(f,g)] \to [K,X] \times_{[K,A]} [K,Y]
	\]
	
	is a surjection. Suppose that $f$ is a fibration, so that every element of $[K,X] \times_{[K,A]} [K,Y]$ can be represented by a pair of maps $u: K \to X, v: K \to Y$ such that $fu = gv := w$. Then the preimage of $(u,v)$ is isomorphic to the set of orbits of $[K \wedge I_+, w]$ under the right action of the group $[K \wedge I_+, u] \times [K \wedge I_+, v]$. In particular, the map is injective iff each of the functions $[K \wedge I_+, u] \times [K \wedge I_+,v] \to [K \wedge I_+, w]$ is surjective.
	
\end{Lemma}

\begin{proof}
	This follows from the same arguments as in \cite[Proposition 2.2.2]{MP12}, where the result is proved in the special case when $u$ and $v$ are nullhomotopic.
\end{proof} 

Next, and along the same lines as \cite[V.5.1,VI.7.1]{BK72}, we seek to understand how the groups $[K \wedge I_+, f]$ behave with respect to $\mathbf{T}$-localisation.

\begin{Lemma} \label{name90} Let $\mathbf{T}$ be a constant localisation system. Let $K$ be a finite based $G$-CW complex, let $X$ be a nilpotent $G$-space, and let $f: K \to X$ be a map. Then:
	
	i) $[K \wedge I_+, f]$ is a nilpotent group, which is finitely $T$-generated (see \cite[Definition 5.6.3]{MP12}) if, for every $i \geq 2$ and $H$, $\pi_i(X^H)$ is finitely $T$-generated, \\
	ii) if the coefficient of $\mathbf{T}$ is 0, then $[K \wedge I_+, f] \to [K \wedge I_+, \phi_{\mathbf{T}} f]$ is $T$-localisation, where $\phi_{\mathbf{T}}$ is a $\mathbf{T}$-localisation of $X$, \\
	iii) if the coefficient of $\mathbf{T}$ is 1, $\mathbb{H}_{\mathbf{T}} \underline{\pi}_1 (X) = 0$, and, for every $i \geq 2$ and $H$, $\pi_i(X^H)$ is finitely $T$-generated, then $[K \wedge I_+, f] \to [K \wedge I_+, \phi_{\mathbf{T}} f]$ is $T$-completion. \\
iv) if, for every $i \geq 2$ and $H$, $\pi_i(X^H)$ is an $f \widehat{\mathbb{Z}}_T$-nilpotent group (resp. $T$-complete), then $[K \wedge I_+, f]$ is $f \widehat{\mathbb{Z}}_T$-nilpotent (resp. $T$-complete).

 \end{Lemma}

\begin{proof}
	This follows by induction up the CW structure on $K$, using Lemma \ref{name95}. In more detail, part i) follows from \cite[Lemma 3.1.3]{MP12} and the fact that a nilpotent group $G$ is finitely $T$-generated iff $G_T$ is $f\mathbb{Z}_T$-nilpotent. Part ii) follows from \cite[Corollary 5.4.11]{MP12}. Part iii) follows from \cite[Corollary 10.4.5]{MP12}, and the condition that $\mathbb{H}_{\mathbf{T}} \underline{\pi}_1 (X) = 0$ ensures that $\underline{\pi}_2(X) \to \underline{\pi}_2(X_{\mathbf{T}})$ is $\mathbf{T}$-localisation, by Theorem \ref{name12}. Part iv) follows from \cite[Lemmas 2.1,2.3]{R24}. 
\end{proof}

\begin{Example} \label{name194}
	If $V$ is a finite dimensional real representation containing at least one copy of the trivial representation, then $S^V$ is a $G$-connected finite $G$-CW complex (see, for instance, ncatlab's page on representation spheres). As we will see in Lemma \ref{name195}, this implies that $S^V$ is based $G$-homotopy equivalent to a finite \textit{based} $G$-CW complex. Therefore, Lemma \ref{name90} implies, for example, that if $X$ is a nilpotent $G$-space and $\mathbf{T}$ is a constant localisation system taking value $(T,0)$, then $\pi_{V+1} (X) \to \pi_{V+1}(X_{\mathbf{T}})$ is $T$-localisation.
\end{Example}

\begin{Lemma} \label{name195}
	Let $X$ be a based $G$-space which is also a $G$-connected finite $G$-CW complex. Then $X$ is based $G$-homotopy equivalent to a finite \textit{based} $G$-CW complex. 
\end{Lemma}

\begin{proof}
	By standard arguments involving the well-pointedness of $G$-CW complexes, it suffices to construct an unbased $G$-homotopy equivalence from $X$ to a finite based $G$-CW complex. Since any basepoint of $X$ is in the interior of some fixed cell, we may as well arrange for the basepoint to be a $0$-cell. Suppose that $X$ has more than one $0$-cell. Let $H$ be a maximal subgroup, in the sense $H$ is not strictly subconjugate to another subgroup satisfying the same property, such that there is a $1$-cell attached along a map from $G/H \times S^0$, with one boundary the basepoint and the other boundary not the basepoint. Since $X$ is finite and $G$-connected, such an $H$ exists. We claim that the induced map $C(G/H) \to X$, from the cone on $G/H$, is injective. If not, then two points on the boundary $G/H$ are sent to the same element $x \in X$. Then the isotropy group of $x$, $K$, is strictly larger than $H$. However, $X^K$ is connected and $(G/H)^K = \emptyset$, so there must be another $1$-cell with boundary $*$ with isotropy group is larger than $H$, a contradiction. Next, consider the pushout:
	
\[	\begin{tikzcd}
		C(G/H) \arrow[tail]{r} \arrow[swap]{d}{\simeq} & X \arrow{d} \\
		* \arrow{r} & X^{'}
	\end{tikzcd} \]

Since the $h$-model structure on $G$-spaces is proper, \cite[Theorem 6.1.2]{R23}, we see that $X \to X^{'}$ is a homotopy equivalence and $X^{'}$ is a finite $G$-CW complex with fewer $0$-cells. In this way, we can reduce to the case where $X$ has a single $0$-cell, the basepoint.

Net we form a `thickened' $G$-CW complex, $Y$, which includes $X$ as a subcomplex, by forming pushouts of the form:

\[
\begin{tikzcd}
	\sqcup G/H \times S^{k-1} \times I \arrow{r} \arrow{d} & Y^{(k-1)} \arrow{d} \\
	\sqcup G/H \times D^k \times I \arrow{r} & Y^{(k)} 
\end{tikzcd}
\]

for $k \geq 2$, letting $Y^{(1)} = X^{(1)}$.

We define the attaching maps $\alpha: G/H \times S^{k-1} \times I \to Y^{(k-1)}$ inductively as follows. We will require $\alpha$ to restrict to the corresponding attaching map, $\alpha_0$, of $X$ at $0 \in I$. Therefore, the inclusion $0 \to I$ will yield a homotopy equivalence $X \to Y$. We define $\alpha$ so that the restriction, $\alpha_1$, induced by $1 \in I$ defines an attaching map of a based $G$-CW complex, $Z$. Inductively, we have that $Z^{(k-1)} \to Y^{(k-1)}$ is a homotopy equivalence. So $\alpha_0$ is homotopic to a map into $Z^{(k-1)}$, which we can take to be a based map since $Z^{(k-1)}$ is $G$-connected and the basepoints are nondegenerate. Therefore, we can inductively define $Y$ and $Z$, and we have a homotopy equivalence from $X$ to $Z$, as desired.
\end{proof}

Before moving on to derive fracture theorems for homotopy classes, we record, in Lemma \ref{name193}, a not-so-obvious property of the groups $[K \wedge I_+,f]$. Recall that a nilpotent group is said to be $T$-completable if $\mathbb{H}_T G = 0$. It is worth bearing in mind the following counterexample: if $\mathbb{Z}_p \oplus \mathbb{Z}_p \to \widehat{\mathbb{Z}}_p$ is an injection, then its cokernel is not $p$-completable. Note, however, that this map does not factor through $\mathbb{Z}_p$. With this in mind, we begin with some preparatory lemmas:

\begin{Lemma} \label{name190}
	If $G$ is a $T$-local nilpotent group and the image of $G$ in $\mathbb{E}_T G$ is normal, then $\mathbb{E}_T G / im(G)$ is a rational nilpotent group.
\end{Lemma}

\begin{proof}
	This follows from the fracture square \ref{dia1} for $K(G,1)$, since if $F$ is the homotopy fibre of $K(G,1) \to \widehat{K(G,1)}_T$, then the cokernel of $G \to \mathbb{E}_T G$ can be identified with $\pi_0 (F)$, using the results of \cite[Lemma 1.4.7]{MP12}. Now, $\pi_0 (F)$ can also be identified with the cokernel of $G_0 \to (\mathbb{E}_T G)_0$, which is rational, and the naturality of the identification implies that the group structures on $\pi_0 (F)$ agree.
\end{proof}

The next lemma is a straightforward consequence of \cite[Lemma 10.4.4]{MP12}:

\begin{Lemma} \label{name191}
	Suppose that $K \xrightarrow{} G \xrightarrow{q} H$ is an exact sequence of nilpotent groups such that the image of $q$ is normal and $K,G,H$ and $coker(q)$ are $T$-completable. Then $\mathbb{E}_T K \to \mathbb{E}_T G \to \mathbb{E}_T H$ is also exact.
\end{Lemma}

The next lemma allows us to recognise completable cokernels:

\begin{Lemma} \label{name192}
	Suppose that $G \to B$ is a map of $T$-local nilpotent groups, let $A := Ab(G)$ and let $K$ be the kernel of $\mathbb{E}_T G \to \mathbb{E}_T B$. Suppose also that:
	
	i) $B$ is abelian, \\
	ii) $B \to \mathbb{E}_T B$ is injective, \\
	iii) the image of $K$ in $\mathbb{E}_T A / im(A)$ is $T$-divisible ($p$-divisible for every $p \in T$).
	
Then $\mathbb{H}_T Q = 0$, where $Q$ is the cokernel of $G \to B$.
\end{Lemma}

\begin{proof}
	Using $ii)$ and the exact sequence of \cite[Lemma 10.4.4]{MP12}, it is straightforward to reduce to the case where $B$ is already $T$-complete, which from now on we assume. Recall that $T$-complete nilpotent groups are $T$-completable, \cite[Lemma 10.4.6]{MP12}, and that the category of $T$-complete nilpotent groups satisfies the closure properties of \cite[Corollary 2.5]{R24}. Therefore, using \cite[Lemma 10.4.4]{MP12} again shows that $\mathbb{H}_T Q = \mathbb{H}_T (im(\mathbb{E}_T G)/ im(G)) = \mathbb{H}_T (im(\mathbb{E}_T A)/ im(A))$, where the images are in $B$. Note that $im(\mathbb{E}_T A)/ im(A)$ is the cokernel of $K \to \mathbb{E}_T A / im(A) := C$. By Lemma \ref{name190}, $C$ is rational, so $\mathbb{H}_T C = \mathbb{E}_T C = 0$ by \cite[Proposition 10.4.7]{MP12}. Therefore, condition iii) and \cite[Proposition 10.4.7]{MP12} imply that $\mathbb{H}_T Q = 0$. 
\end{proof}

We can now prove the property of the groups $[K \wedge I_+, f]$ that we've been aiming for:

\begin{Lemma} \label{name193}
	Let $T$ be a set of primes, and let $\mathbf{T}$ denote the constant localisation system taking value $(T,1)$. Suppose that $\mathbf{L}$ is a localisation system such that the underlying set of primes of $\mathbf{L}$ is constantly $T$. Let $f: K \to X$ be a map, where $X$ is a nilpotent $G$-space and $K$ is a finite based $G$-CW complex. Suppose that $\mathbb{H}_{\mathbf{T}} \underline{\pi}_1 (X) = 0$, and, for every $i \geq 2$ and $H$, $\pi_i(X^H)$ is finitely $T$-generated. Then:
	
	i) $[K \wedge I_+, \phi_{\mathbf{L}} f]$ is $T$-completable, \\
	ii) $[K \wedge I_+, \phi_{\mathbf{L}} f] \to [K \wedge I_+, \phi_{\mathbf{T}} f]$ is $T$-completion.
\end{Lemma}

\begin{proof}
	Due to Lemma \ref{name90} ii), we may as well assume that $X$ is $T$-local. We will induct on the cells of $K$, where recall that $K$ has a single $0$-cell which is the basepoint. Suppose that the result is true for a subcomplex $E$ and we attach a cell to $E$ via a cofibre sequence $S \to E \to L$, where $S$ is of the form $(G / H)_+ \wedge S^n$ with $n \geq 1$ (or $n = 0$ and we take $S$ to be a finite wedge of such $G$-spaces corresponding to the $1$-cells). Consider the diagram of $T$-local nilpotent groups:
	
\[	\begin{tikzcd}
	{[}E \wedge I^2_+, f{]} \arrow{r} \arrow{d} & {[}\Sigma^2 S, X{]} \arrow{r} \arrow{d} &  {[}L \wedge I_+, f{]} \arrow{r} \arrow{d} & {[}E \wedge I_+, f{]} \arrow{r} \arrow{d} & {[}\Sigma S, X{]} \arrow{d} \\
	{[}E \wedge I^2_+, \phi_{\mathbf{L}} f{]} \arrow{r} \arrow{d} & {[}\Sigma^2 S, X_{\mathbf{L}}{]} \arrow{r} \arrow{d} &  {[}L \wedge I_+, \phi_{\mathbf{L}} f{]} \arrow{r} \arrow{d} & {[}E \wedge I_+, \phi_{\mathbf{L}}f{]} \arrow{r} \arrow{d} & {[}\Sigma S, X_{\mathbf{L}}{]} \arrow{d} \\
	{[}E \wedge I^2_+, \phi_{\mathbf{T}} f{]} \arrow{r}  & {[}\Sigma^2 S, X_{\mathbf{T}}{]} \arrow{r}  &  {[}L \wedge I_+, \phi_{\mathbf{T}} f{]} \arrow{r}  & {[}E \wedge I_+, \phi_{\mathbf{T}}f{]} \arrow{r}  & {[}\Sigma S, X_{\mathbf{T}}{]} 
\end{tikzcd} \]
	
	Lemma \ref{name95} implies that all three rows are exact, and Lemma \ref{name90}iii) implies that the composite vertical maps are all $T$-completions. Inductively, all vertical maps, except the middle one, between the second and third row are $T$-completions. By Lemma \ref{name191}, to prove the lemma it suffices to prove that the cokernels of the maps $[E \wedge I^2_+, \phi_{\mathbf{L}} f] \to [\Sigma^2 S, X_{\mathbf{L}}]$ and $[E \wedge I_+, \phi_{\mathbf{L}} f] \to [\Sigma S, X_{\mathbf{L}}]$ are $T$-completable. The proofs are the same so we concentrate on the latter case. The map $[\Sigma S, X_{\mathbf{L}}] \to [\Sigma S, X_{\mathbf{T}}]$ is either an isomorphism or the $T$-completion of an $f \mathbb{Z}_T$-nilpotent group, so is injective. So conditions $i)$ and $ii)$ of Lemma \ref{name192} are satisfied. For condition $iii)$, note that the image of ``K" factors through the cokernel of $Ab([L \wedge I_+, f]) \to \mathbb{E}_T Ab([L \wedge I_+, f]) $, which is rational by Lemma \ref{name190}. So the conditions of Lemma \ref{name192} are satisfied, and so the cokernels are completable, as desired.
\end{proof}

We now prove our main fracture theorems for homotopy classes of maps, starting with:

\begin{Theorem} \label{name91} Let $\mathbf{T}, \mathbf{S}$, and, for each $i$ in some indexing set $I$, $\mathbf{T_i}$ be constant localisation systems such that $\mathbf{T}$ and $\mathbf{S}$ have coefficient 0, $T = \cup_i T_i$ and $T_i \cap T_j = S$, for all $i \neq j$. Let $K$ be a finite based $G$-CW complex and let $X$ be a $T$-local nilpotent $G$-space such that, if $\mathbf{T_i}$ has coefficient 1, then for every subgroup $H$ of $G$, $\mathbb{H}_{T_i}\pi_1 (X^H) = 0$, and for each $i \geq 2$, $\pi_i(X^H)$ is finitely $T_i$-generated. Then the following diagram is a pullback of sets:
	
	\[\begin{tikzcd}
		{[K,X]} \arrow{r} \arrow{d} & {[K,\prod_i X_{\mathbf{T_i}}]} \arrow{d}  \\
		{[K, X_{\mathbf{S}}]} \arrow{r} & {[K, (\prod_i X_{\mathbf{T_i}})_{\mathbf{S}}]}
	\end{tikzcd} \]

\end{Theorem} 

\begin{proof}
	The map from $[K,X]$ to the pullback is surjective by Theorem \ref{name88} and Lemma \ref{name89}. The map is injective by Lemma \ref{name89} and Lemma \ref{name90}. In more detail, to see that the surjectivity hypothesis in the final sentence of Lemma \ref{name89} is satisfied, surjectivity tells us that we can find a map $\mu: K \to X$ projecting onto a representative $(u,v)$ of any element of the pullback. Then we can apply Theorem \ref{name88} to give a fracture square for $K([K \wedge I_+, \mu], 1)$, and the fact that $K([K \wedge I_+, \mu],1)$ is connected tells us that the required map is surjective, via use of Lemma \ref{name90}.
\end{proof}

Note that the fracture theorems for nilpotent groups given in \cite[Theorem 7.2.1 ii), Theorem 12.3.2]{MP12}, are both consequences of Theorem \ref{name91}. If the homotopy groups of $X^H$ are finitely $T$-generated, we only require the underlying set of primes of the $\mathbf{T}_i$ to be constant:

\begin{Theorem} \label{name188}
	Let $T$ be a set of primes and, for every subgroup $H$ of $G$, let $S_H$ be a set of primes. Consider a commutative square of nilpotent $G$-spaces:
	\begin{center}

		\begin{tikzcd}
			X \arrow{r}{f} \arrow[swap]{d}{\phi} & Y \arrow{d}{\psi} \\
			A \arrow[swap]{r}{g} & B \\
		\end{tikzcd}
	\end{center}	
	
	such that the following conditions are satisfied:
	
	i) $X^H,Y^H$ are $T$-local and $A^H,B^H$ are $S_H$-local, \\
	ii) $f^H$ is an $\mathbb{F}_{T}$-equivalence and $g^H$ is an $\mathbb{F}_{S_H}$-equivalence, \\
	iii) $\phi^{H}, \psi^{H}$ are $\mathbb{Q}$-equivalences, \\
	iv) for every $H$, $\mathbb{H}_T \pi_1(X^H) = 0$ and $\pi_i(X^H)$ is finitely $T$-generated for $i \geq 2$.
	
	Then, for every finite based $G$-CW complex $K$, there is a pullback of sets:
		\[		
	\begin{tikzcd}
		{[K,X]} \arrow{r}{f_*} \arrow[swap]{d}{\phi_*} & {[K,Y]} \arrow{d}{\psi_*} \\
		{[K,A]} \arrow[swap]{r}{g_*} & {[K,B]} \\
	\end{tikzcd}
	\]
\end{Theorem}

\begin{proof}
	Suppose that we have a diagram of sets:
	
\[	\begin{tikzcd}
		R \arrow{r} \arrow{d} & S \arrow{r} \arrow{d} & L \arrow{d} \\
		U \arrow{r} & V \arrow{r} & W
	\end{tikzcd} \]

and let $P_1$ denote the left hand pullback, $P_2$ denote the right hand pullback and $P$ denote the composite pullback. The key point is that if  $R \to P$ is an isomorphism, then $R \to P_1$ is injective. Therefore, we can use a pasting argument to prove the theorem, since in the situation at hand the map $R \to P_1$ will also be surjective via Theorem \ref{name88} and Lemma \ref{name89}. The starting point for the proof is the case of Theorem \ref{name91} where $I = \{i\}$, $\mathbf{T}_i$ takes value $(T,1)$ and $\mathbf{S}$ is rationalisation.
\end{proof}

\begin{Example}
	
To finish this subsection, we give an example to show that the square:  

\[\begin{tikzcd}
	{[K,X]} \arrow{r} \arrow{d} & {[K,X_{\mathbf{T}}]} \arrow{d}  \\
	{[K, X_0]} \arrow{r} & {[K, (X_{\mathbf{T}})_0]}
\end{tikzcd} \]

need not be a pullback of sets if the underlying set of primes of $\mathbf{T}$ is not constant, where $K$ is a finite based $G$-CW complex and, for every $H$, $X^H$ is $f \mathbb{Z}_{T(\orb{H})}$-nilpotent. Of course, the comparison map $[K,X] \to [K,X_{\mathbf{T}}] \times_{[K, (X_{\mathbf{T}})_0]} [K,X_0]$ is always a surjection by Theorem \ref{name88} and Lemma \ref{name89}. We let $G = C_2$, and let $\mathbf{T}(\orb{G}) = (\{p,q\},1)$, and $\mathbf{T}(\orb{e}) = (\{p\},1)$. We let $X = S^{n+2}_{\mathbf{L}}$, where $\mathbf{L}$ denotes the localisation system induced by $\mathbf{T}$ taking values with coefficient $0$, so $\pi_{n+2}(X^G) = \mathbb{Z}_{\{p,q\}}$ and $\pi_{n+2}(X^e) = \mathbb{Z}_{\{p\}}$. We let $K$ be the cofibre:
\[ \Sigma S^n \wedge (G_+) \to \Sigma S^n \wedge (e_+) \to K
\]
where the first map is induced by the constant map $G \to e$. Suppose also that $\pi_{n+3}(S^{n+2})$ is a finite group with no $p$-torsion. The cofibre sequence implies that there is a map of short exact sequences:

\[\begin{tikzcd}
	0 \arrow{r} & {[\Sigma K, X_{\mathbf{T}}]} \oplus {[\Sigma K, X_0]} \arrow{r} \arrow{d} & \widehat{\mathbb{Z}}_{\{p,q\}} \oplus \mathbb{Q} \arrow[two heads]{d} \arrow{r} & \widehat{\mathbb{Z}}_p \oplus \mathbb{Q} \arrow[two heads]{d} \arrow{r} & 0 \\
	0 \arrow{r} & {[\Sigma K, (X_{\mathbf{T}})_0]} \arrow{r} & \widehat{\mathbb{Q}}_{\{p,q\}} \arrow{r} & \widehat{\mathbb{Q}}_p \arrow{r} & 0 \\
\end{tikzcd}
\]

It follows that the map $[\Sigma K, X_{\mathbf{T}}] \oplus [\Sigma K, X_0] \to [\Sigma K, (X_{\mathbf{T}})_0]$ can be identified with the rationalisation $\widehat{\mathbb{Z}}_q \to \widehat{\mathbb{Q}}_q$, this being the map between the kernels above since $\widehat{\mathbb{Z}}_{\{p,q\}} \cong \widehat{\mathbb{Z}}_p \oplus \widehat{\mathbb{Z}}_q$, which is not surjective. It follows that the square above is not a pullback of sets. In particular, there exist maps $f,g: K \to X$, such that $f_{\mathbf{T}} \simeq g_{\mathbf{T}}$ and $f_0 \simeq g_0$, but $f$ and $g$ are not homotopic. Note also that if $H$ is a subgroup of $G$, then $f^H_{\mathbf{T}(\orb{H})} \simeq g^H_{\mathbf{T}(\orb{H})}$, $f^H_0 \simeq g^H_0$, $K^H$ is a finite based CW-complex, and $X^H$ is $f \mathbb{Z}_{T(\orb{H})}$-nilpotent, so $f^H \simeq g^H$ by Theorem \ref{name91}. Therefore, we have constructed two equivariant maps which are non-equivariantly homotopic after taking fixed points at any subgroup, but which are not equivariantly homotopic. In fact, by inspection of Lemma \ref{name89} and what we have just proved, we may as well take $g$ to be equivariantly nullhomotopic.

\end{Example}

\subsection{Equivariant phantom maps} \label{name186}

In this subsection, we show that the classical application of completion to the study of phantom maps with nilpotent codomain, \cite[Lemma 11.6.1]{MP12}, continues to hold in the current equivariant setting. Equivariant phantom maps have already been studied in \cite{S88} and \cite{P01}, so our discussion will be brief. The key point is the following algebraic lemma:

\begin{Lemma} \label{name171}
	Let $... \to G_2 \to G_1 \to G_0$ be a sequence of nilpotent groups and $T$ be a set of primes. Then $lim^1 G_i = *$ if either of the following conditions is satisfied: \\
	i) every $G_i$ is $f \mathbb{Q}$-nilpotent, \\
	ii) every $G_i$ is $f \widehat{\mathbb{Z}}_T$-nilpotent.
\end{Lemma}

\begin{proof}
	i) This follows from \cite[Lemma 6.8.2]{MP12}. \\
	ii) This is a well-known result, see \cite[Theorem 2.6]{S77} and \cite[Lemma 11.6.1]{MP12}, but we take the opportunity to present an alternative argument. First suppose that all of the $G_i$ are abelian. Then $lim^1 G_i$ is $T$-complete, since it is the cokernel of a map between $T$-complete abelian groups, see \cite[pg. 148]{M99} and \cite[Lemma 2.3]{R24}. Moreover, if $p \in T$ then $G_i / pG_i$ is a $T$-torsion $f \widehat{\mathbb{Z}}_T$-nilpotent group, and so is finite by \cite[Lemma 2.6]{R24}. Therefore, the Mittag-Leffler condition and the exact sequence of \cite[Proposition 2.3.1]{MP12} imply that $lim^1 G_i/ pG_i = 0$ and $lim^1 G_i \xrightarrow{.p} \lim^1 G_i$ is surjective. So $lim^1 G_i$ is $p$-divisible for every $p \in T$, and since it is also $T$-complete it must vanish by \cite[Proposition 10.4.7]{MP12}. 
	
	Now recall that a $T$-complete nilpotent group, $G$, has a functorial $T$-complete central series, $\{\Lambda^i G\}$, where $\Lambda^i G$ is the image of $\widehat{(\Gamma_i G)}_T \to G$, for $\Gamma_i G$ the ith term of the lower central series of $G$. Use of this central series and the exact sequence of \cite[Proposition 2.3.1]{MP12} shows that $lim^1 G_i = *$ whenever there is a common bound on the nilpotency classes of the $G_i$.
	
	For the more general case, consider the homotopy limit of the doubly indexed diagram of spaces $\{ K(G_i / \Lambda^j G_i,1)\}_{i,j}$, which is easily seen to coincide with the homotopy limit of the tower $\{K(G_i,1)\}$, see \cite[Example 4.3]{BK72}. We want to show that $lim^1 G_i = *$ which is equivalent to this homotopy limit being connected by \cite[Proposition 2.2.9]{MP12}. We've already shown that $holim_i G_i / \Lambda^j G_i$ is connected for every $j$. Moreover, the exact sequence $1 \to \Lambda^j G_i / \Lambda^{j+1} G_i \to G_i / \Lambda^{j+1} G_i \to G_i / \Lambda^{j} G_i \to 1$ induces an exact sequence on limits which shows that $lim_i G_i / \Lambda^{j+1} G_i \to lim_i G_i / \Lambda^j G_i$ is surjective. Therefore, the Mittag-Leffler condition implies that the homotopy limit of the doubly indexed diagram is connected, as desired.
\end{proof}

To apply Lemma \ref{name171} to the study of phantom maps, we let $K$ be a \textit{based} $G$-CW-complex with finitely many cells of each dimension. Then, using the methods of \cite[Proposition 2.1.9]{MP12}, for any $G$-space $X$ there is a surjection $[K,X] \to lim[K^{(n)}, X]$ where the premiage of $(u_n) \in lim[K^{(n)},X]$ can be identified with $lim^1 [K^{(n)} \wedge I_+, u_n]$. Therefore, Lemmas \ref{name95} and \ref{name171} imply the following theorem:

\begin{Theorem} \label{name172}
	Let $K$ be a based $G$-CW complex with finitely many cells in each dimension (ie of finite type) and let $T$ be a set of primes. If $X$ is a $G$-space such that either:
	
	 i) $\pi_i(X^H)$ is  $f \mathbb{Q}$-nilpotent for all $i \geq 2$ and $H$, or \\ 
	 ii) $\pi_i(X^H)$ is $f \widehat{\mathbb{Z}}_T$-nilpotent for all $i \geq 2$ and $H$,
	 
	then $[K,X] \to lim[K^{(n)}, X]$ is a bijection. Therefore, if $f: K \to X$ is phantom (ie $f|_{K^{(n)}}$ is nullhomotopic for all $n$), then $f$ is nullhomotopic.
\end{Theorem}

Note that it makes sense for the above result to only depend on $\pi_i(X^H)$ for $i \geq 2$ since, as observed on \cite[pg. 135]{Z87}, non-equivariantly a map $f: K \to X$ is phantom iff $f$ factors through the universal cover of $X$ and $\tilde{f}: K \to \tilde{X}$ is phantom. We now assume that $X$ is a nilpotent $G$-space such that $X^H$ is $f \mathbb{Z}_T$-nilpotent for all $H$. As in \cite[Theorems 8.2.1/13.1.1]{MP12}, it can be shown that if $L$ is a finite based $G$-CW complex, then $[L,X] \to [L, \widehat{X}_T]$ is an injection. Here, $\widehat{X}_T$ denotes $X_{\mathbf{T}}$, for $\mathbf{T}$ the constant localisation system with underlying set of primes $T$ and coefficient $1$. Note that the `neighbourhood groups' which appear in the proof arise from extending the exact sequence of Lemma \ref{name95} to include $\pi_0$ terms, as in \cite[V.5]{BK72}, and so share similar properties to those of Lemma \ref{name90}. It follows that if $K$ is a based $G$-CW complex of finite type, then a map $f: K \to X$ is phantom iff the composite $\widehat{\phi} f: K \to \widehat{X}_T$ is nullhomotopic.

\subsection{Nilpotent $G$-spaces and Postnikov towers} \label{name80}

We now discuss the relationship between $\mathbf{T}$-localisations of nilpotent $G$-spaces and towers of principal fibrations. The arguments of this section are similar to those of \cite{M82}, except that we derive the results appropriate to non-constant localisation systems, and make a distinction between nilpotent $G$-spaces and bounded nilpotent $G$-spaces, and the types of Postnikov tower they are equivalent to. First, we will define what it means for a $\underline{\pi}$-group to be $\mathcal{B}$-nilpotent, where $\mathcal{B}$ is a class of coefficient systems. Then we will define the analogue of Postnikov towers in the equivariant setting, and we will show that a $G$-space, X, is equivalent to a \textit{weak} Postnikov $\mathcal{B}$-tower iff its homotopy groups are $\mathcal{B}$-nilpotent $\underline{\pi}$-groups, where $\underline{\pi} = \underline{\pi}_1(X)$. Finally, we will show that $\mathbf{T}$-local nilpotent $G$-spaces are equivalent to weak Postnikov $\mathcal{B}_{\mathbf{T}}$-towers, where $\mathcal{B}_{\mathbf{T}}$ is the class of $\mathbf{T}$-local coefficient systems.

\begin{Definition}  \label{name24}
	Let $\mathcal{B}$ be a class of coefficient systems of abelian groups. Let $\underline{\pi}$ be a coefficient system of groups and let $\mathbf{G}$ be a coefficient system of groups admitting an action of $\underline{\pi}$ by automorphisms. We say that $\mathbf{G}$ is a $\mathcal{B}$-nilpotent $\underline{\pi}$-group if there is a descending sequence of normal $\underline{\pi}$-subgroups:
	
	\begin{center}
		$\mathbf{G} = \mathbf{G_0} \supseteq \mathbf{G_1} \supseteq \mathbf{G_2} \supseteq ... $
	\end{center}
	
	such that:
	
	i) $\underline{\pi}$ acts trivially on $\frac{\mathbf{G_{i-1}}}{\mathbf{G_i}}$, \\
	ii) $\frac{\mathbf{G_{i-1}}}{\mathbf{G_i}} \in \mathcal{B}$, \\
	iii) for every $H$, $\frac{\mathbf{G_{i-1}}}{\mathbf{G_i}}(\orb{H}) \to \frac{\mathbf{G}}{\mathbf{G_i}}(\orb{H})$ has central image,\\
	iv) for every $H$, $\frac{\mathbf{G_{i-1}}}{\mathbf{G_i}}(\orb{H}) = 0$ for all but finitely many $i$. 
	
\end{Definition}

\begin{Definition} Call a $\mathcal{B}$-nilpotent $\underline{\pi}$-group \textit{bounded} if the filtration in Definition \ref{name24} can be replaced by a finite filtration terminating at $1$. \end{Definition} 

\begin{Definition}
	We call a $G$-space $X$ $\mathcal{B}$-nilpotent if it is $G$-connected and, for all $i \geq 1$, $\underline{\pi}_i(X)$ is a $\mathcal{B}$-nilpotent $\underline{\pi}_1(X)$-group. We say that a $\mathcal{B}$-nilpotent $G$-space, $X$, is \textit{bounded} if the homotopy groups $\underline{\pi}_i(X)$ are all bounded $\mathcal{B}$-nilpotent $\underline{\pi}_1(X)$-groups.
\end{Definition}

Note that a $G$-space is nilpotent iff it is $\mathcal{A}$-nilpotent, where $\mathcal{A}$ is the class of all coefficient systems of abelian groups. This follows from the fact that if $X$ is a nilpotent space, then there are functorial filtrations of $\pi_i(X)$ satisying the conditions of the previous definition - the lower central series when $i = 1$, and the filtration induced by the augmentation ideal, $\{I^n \pi_i(X)\}$, for $i \geq 2$.

\begin{Definition}
	A map of $G$-spaces, $f:X \to  Y$, is called a principal $K(\underline{A},n)$-fibration if it is the pullback of the path-space fibration along a map $k: Y \to K(\underline{A},n+1)$. In particular, $f$ is a fibration with fibre $K(\underline{A},n)$.
\end{Definition}

\begin{Definition} \label{name27}
	Let $\mathbf{Q}$ be the totally ordered set consisting of pairs of natural numbers ordered by $(m,n) \leq (p,q)$ iff $m < p$ or $m = p$ and $n \leq q$. A weak Postnikov $\mathcal{B}$-tower is a functor $\mathbf{Q} \to G \text{--} Sp$, where $G \text{--} Sp$ is the category of $G$-spaces, satisfying:
	
	i) $X_{1,1} = *$, \\
	ii) $X_{n+1,1} \to \lim_i X_{n,i}$ is a weak equivalence, \\
	iii) The map $X_{n,i+1} \to X_{n,i}$ is a principal $K(\underline{B}_{n,i}, n)$-fibration for some $\underline{B}_{n,i} \in \mathcal{B}$, \\
	iv) for every $n$ and $H$, $X_{n,i+1}^H \to X_{n,i}^H$ is a weak equivalence for all but finitely many $i$.
\end{Definition}

\begin{Definition}
	A Postnikov $\mathcal{B}$-tower is a weak Postnikov $\mathcal{B}$-tower such that the maps $X_{n+1,1} \to \lim_i X_{n,i}$ of condition ii) above are $G$-homeomorphisms. 
\end{Definition}

We have the principal fibration lemma:

\begin{Lemma} \label{name28}
	Let $f:X \to Y$ be a map of well-pointed $G$-connected $G$-spaces with the homotopy type of a $G$-CW complex, such that $Ff \simeq K(\underline{A},n)$ for some coefficient system $\underline{A}$ and $n \geq 1$. Then $\underline{\pi}_1(X)$ acts trivially on $\underline{\pi}_*(Ff)$ iff there is a weak equivalence $X \to Fk$ over $Y$, for some cofibration $k: Y \to K(\underline{A},n+1)$. 	
\end{Lemma}

\begin{proof}
	See, for instance, \cite[Lemma 2.3.10]{R23}.
\end{proof}

\begin{Lemma} A $G$-space is $\mathcal{B}$-nilpotent iff it is weakly equivalent to a weak Postnikov $\mathcal{B}$-tower.
\end{Lemma}

\begin{proof}
	To see that weak Postnikov $\mathcal{B}$-towers are $\mathcal{B}$-nilpotent, it suffices to show that $\underline{\pi}_n(X)$ is a $\mathcal{B}$-nilpotent $\underline{\pi}_1(X)$-group. Let $\mathbf{G}_i$ be the kernel of $\underline{\pi}_n(X) \to \underline{\pi}_n(X_{n,i+1})$. The quotients $\frac{\mathbf{G}_{i-1}}{\mathbf{G_{i}}}$ correspond to the coefficient systems $\underline{B}_{n,i}$ appearing in the tower, and so these are in $\mathcal{B}$ by assumption. Now, $\underline{\pi}_1(X)$ acts trivially on $\frac{\mathbf{G}_{i-1}}{\mathbf{G}_i}$ by Lemma \ref{name28}, and the required inclusions are central for the same reason as in the non-equivariant case, namely \cite[Lemma 1.4.7 v)]{MP12}. Finally, $\frac{\mathbf{G_{i-1}}}{\mathbf{G_i}}(\orb{H}) = 0$ for all but finitely many $i$, since $X_{n,i+1}^H \to X_{n,i}^H$ is a weak equivalence for all but finitely many $i$.
	
	Next assume that a $G$-CW complex, $X$, is $\mathcal{B}$-nilpotent. So, for each $n$, we have a filtration of $\underline{\pi}_n(X)$, $\{\mathbf{G}^n_i\}$, satisfying the conditions of Definition \ref{name24}. We define $X^0_{n,i+1}$ by first attaching cells to $X$ along all possible maps $(\frac{G}{H})_+ \wedge S^n \to X$ representing an element of $\mathbf{G}^n_i (\orb{H})$, for some $H \leq G$. Then, inductively define $X^{j}_{n,i+1}$, for each $j \geq 1$, by attaching a cell to $X^{j-1}_{n,i+1}$ along every possible map $(\frac{G}{H})_+ \wedge S^{n+j} \to X^{j}_{n,i+1}$, for any $H \leq G$. Define $X_{n,i+1}$ as the union of the $X_{n,i+1}^j$. Then:
	
	i) $\underline{\pi}_j(X_{n,i+1}) = \underline{\pi}_j (X)$ for $j < n$, \\
	ii) $\underline{\pi}_n(X_{n,i+1}) = \frac{\underline{\pi}_n(X)}{\mathbf{G}_i}$, \\
	iii) $\underline{\pi}_j(X_{n,i+1}) = 0$ for $j > n$.
	
	Moreover, we have an inclusion $X_{n,i} \to X_{m,j}$, whenever $(m,j) \leq (n,i)$ in $\mathbf{Q}$. Since the action of $\underline{\pi}_1(X)$ is trivial on each $\frac{\mathbf{G}_{i-1}}{\mathbf{G}_i}$, each of the maps $X_{n,i+1} \to X_{n,i}$ is equivalent to a $K(\underline{B}_{n,i},n)$-principal fibration, with $\underline{B}_{n,i} \in \mathcal{B}$. To define the weak Postnikov $\mathcal{B}$-tower, we keep each $X_{n,1}$ fixed, and inductively replace each $X_{n,i}$, for $i \geq 2$, using Lemma \ref{name28}. Results of Waner, \cite[Corollary 4.14]{W80}, imply that by doing this we never leave the category of well-pointed $G$-spaces with the homotopy type of a $G$-CW complex. Unfortunately, we could leave this category by taking inverse limits, which is why we leave each $X_{n,1}$ fixed and only require a weak equivalence in Definition \ref{name27} ii).
\end{proof}

If we restrict attention to bounded $\mathcal{B}$-nilpotent $G$-spaces, such as pointwise simply connected $G$-spaces, then the same proof shows:

\begin{Lemma}A bounded $\mathcal{B}$-nilpotent $G$-space is weakly equivalent to a Postnikov $\mathcal{B}$-tower.	
\end{Lemma}

As promised, we have the following characterisation of $\mathbf{T}$-local nilpotent spaces:

\begin{Theorem} \label{name185}
	A nilpotent $G$-space is $\mathbf{T}$-local iff it is $\mathcal{B}_{\mathbf{T}}$-nilpotent, where $\mathcal{B}_{\mathbf{T}}$ is the class of $\mathbf{T}$-local coefficient systems. A $G$-space is a $\mathbf{T}$-local bounded $\mathcal{A}$-nilpotent $G$-space iff it is bounded $\mathcal{B}_{\mathbf{T}}$-nilpotent. 
\end{Theorem}

\begin{proof}
	If $X$ is $\mathcal{B}_{\mathbf{T}}$-nilpotent it is easily verified that all of the homotopy groups of $X^H$ are $\mathbf{T}(\orb{H})$-local, which implies that $X$ is $\mathbf{T}$-local. Suppose that $X$ is $\mathbf{T}$-local. Recall that we have a central $\underline{\pi}_1(X)$-series for $\underline{\pi}_i(X)$ induced by the functorial lower central series when $i = 1$, or the functorial augmentation ideal series when $i \geq 2$. Localising these series at $\mathbf{T}$, and considering their images in $\underline{\pi}_i(X)$, expresses each $\underline{\pi}_i(X)$ as a $\mathcal{B}_{\mathbf{T}}$-nilpotent $\underline{\pi}_1(X)$-series, and so $X$ is $\mathcal{B}_{\mathbf{T}}$-nilpotent. If $X$ is bounded $\mathcal{A}$-nilpotent, then the lower central series terminates after finitely many stages and so $X$ is bounded $\mathcal{B}_{\mathbf{T}}$-nilpotent.
\end{proof}

\subsection{Localisation at equivariant cohomology theories} \label{misc222}

We end the paper by tying up the following loose end. Namely, in \cite[Ch. II]{M96}, localisations of nilpotent $G$-spaces were defined relative to equivariant cohomology theories, and we would like to compare this to our localisations relative to $\mathbf{T}$-equivalences.

\begin{Definition}
	A map of $G$-spaces, $f:X \to Y$, is a cohomology $\mathbf{T}$-equivalence if for all $\mathbf{T}$-local coefficient systems $\underline{A}$, $f^*: H^*(Y;\underline{A}) \to H^*(X; \underline{A})$ is an isomorphism.
\end{Definition}

We have:

\begin{Theorem} \label{name29} A map of $G$-spaces, $f:X \to Y$, is a $\mathbf{T}$-equivalence iff it is a cohomology $\mathbf{T}$-equivalence.	
\end{Theorem}

\begin{proof}
	If $f$ is a $\mathbf{T}$-equivalence, then, since each $K(\underline{A},n)$ is $\mathbf{T}$-local for every $\mathbf{T}$-local coefficient system $\underline{A}$, $f$ is a cohomology $\mathbf{T}$-equivalence. 
	
	If $f$ is a cohomology $\mathbf{T}$-equivalence, then we can assume that $X,Y$ are well-pointed. Now, $\Sigma^2 f: \Sigma^2 X \to \Sigma^2 Y$ is a cohomology $\mathbf{T}$-equivalence between pointwise simply connected $G$-spaces. By the first part, $(\Sigma^2 f)_{\mathbf{T}}$ is also a cohomology $\mathbf{T}$-equivalence between pointwise simply connected $\mathbf{T}$-local $G$-spaces. Simply connected $\mathbf{T}$-local $G$-spaces are weakly equivalent to strict Postnikov $\mathcal{B}_{\mathbf{T}}$-towers by the results of Subsection \ref{name80}, so the equivariant analogue of co-HELP, \cite[Theorem 3.3.7]{MP12}, implies that $(\Sigma^2 f)_{\mathbf{T}}$ is a weak equivalence and, therefore, that $\Sigma^2 f$ is a $\mathbf{T}$-equivalence. It follows that $f$ is a $\mathbf{T}$-equivalence. 	
\end{proof}

\begin{Corollary}
	Localisation with respect to $\mathbf{T}$-equivalences is equivalent to localisation with respect to cohomology $\mathbf{T}$-equivalences.
\end{Corollary}	

\bibliography{References}
\bibliographystyle{alpha}

\end{document}